\documentclass[a4paper,10pt]{scrartcl}
\pdfoutput=1
\usepackage{ifpdf}
\usepackage[english]{babel}
\usepackage[utf8]{inputenc}
\usepackage[T1]{fontenc}
\usepackage{lmodern}
\usepackage{enumerate}
\usepackage{geometry}
\usepackage{amsfonts,amsmath,amssymb,amsopn}
\usepackage{amsthm}
\usepackage{mathtools}
\usepackage{stmaryrd}
\usepackage{nicefrac}
\usepackage{exscale}

\usepackage[numbers,square]{natbib}
\usepackage{url} 
\usepackage[colorlinks=true,pdfpagelabels,unicode]{hyperref}
\hypersetup{linkcolor=blue, urlcolor=blue, citecolor=blue} 

\allowdisplaybreaks 

\theoremstyle{plain}

\newtheoremstyle{theo}
	{3pt} 
	{3pt} 
	{\itshape} 
	{} 
		{\bfseries} 
	{\\} 
	{ } 
	{\thmname{#1}\thmnumber{ #2.}\thmnote{ - #3}} 
\theoremstyle{theo}

\newtheorem{definition}{Definition}[section]
\newtheorem{lemma}[definition]{Lemma}
\newtheorem{theorem}[definition]{Theorem}

\newtheorem{proposition}[definition]{Proposition}

\newenvironment{bew}{\begin{proof}[\bfseries Proof:]}{\end{proof}}

\newtheoremstyle{remark}
	{3pt} 
	{3pt} 
	{} 
	{} 
		{\bfseries} 
	{} 
	{ } 
	{\thmname{#1}\thmnumber{ #2.}\thmnote{ - #3}} 
\theoremstyle{remark}

\DeclareMathOperator{\bomega}{\overline{\Omega}}
\DeclareMathOperator{\romega}{\partial\Omega}

\DeclareMathOperator{\intd}{d\!}

\DeclareMathOperator{\wto}{\rightharpoonup}
\DeclareMathOperator{\wsto}{\stackrel{\star}{\wto}}
\newcommand{\fracpow}{\varrho}
\newcommand{\epsi}{\varepsilon}

\newcommand{\DA}[1][\fracpow]{D\!\left(A^{#1}\right)} 

\newcommand{\HP}{\mathcal{P}}

\newcommand{\Tme}{T_{max,\;\!\epsi}}
\newcommand{\GNI}{Gagliardo--Nirenberg inequality}

\newcommand{\into}[1]{\int_0^{#1}\!}
\newcommand{\intot}{\into{t}}
\newcommand{\intoT}{\into{T}}

\newcommand{\intomega}{\int_{\Omega}\!} 

\newcommand{\inttminusomega}{\int_{(t-1)_+}^t\!\intomega}
\newcommand{\intoTomega}{\intoT\!\intomega}

\newcommand{\intinfomega}{\int_0^\infty\!\!\intomega}

\newcommand{\intromega}{\int_{\romega}\!} 

\newcommand{\Lo}[1][1]{L^{#1}(\Omega)} 
\newcommand{\W}[1][1,2]{W^{#1}(\Omega)}

\newcommand{\LSp}[2]{L^{#1\;\!}\!\left(#2\right)} 
\newcommand{\LSpn}[2]{L^{#1\;\!}\!(#2)}
\newcommand{\LSpb}[2]{L^{#1\;\!}\!\big(#2\big)}
\newcommand{\LSploc}[2]{L_{loc}^{#1}\!\left(#2\right)} 

\newcommand{\LSplocb}[2]{L_{loc}^{#1}\big(#2\big)}

\newcommand{\CSp}[2]{C^{#1}\!\left(#2\right)}

\newcommand{\R}{\mathbb{R}}
\newcommand{\N}{\mathbb{N}}

\newcommand{\uep}{u_\epsi}
\newcommand{\nep}{n_\epsi}
\newcommand{\cep}{c_\epsi}
\newcommand{\fep}{f_\epsi}
\newcommand{\gep}{g_\epsi}
\newcommand{\rhoep}{\rho_\epsi}
\newcommand{\cs}{c_*}
\newcommand{\hcep}{\hat{c}_{\epsi}}
\newcommand{\hc}{\hat{c}}
\newcommand{\dimN}{N}

\makeatletter
\def\@fnsymbol#1{\ensuremath{\ifcase#1\or *\or \ddagger\or \#\or
   \mathsection\or \mathparagraph\or \|\or **\or \dagger\dagger
   \or \ddagger\ddagger \else\@ctrerr\fi}}
   
 \def\@lbibitem[#1]#2#3{%
  \if\relax\@extra@b@citeb\relax\else
    \@ifundefined{br@#2\@extra@b@citeb}{}{%
     \@namedef{br@#2}{\@nameuse{br@#2\@extra@b@citeb}}%
    }%
  \fi
  \@ifundefined{b@#2\@extra@b@citeb}{%
   \def\NAT@num{}%
  }{%
   \NAT@parse{#2}%
  }%
  \def\NAT@tmp{#1}%
  \expandafter\let\expandafter\bibitemOpen\csname NAT@b@open@#2\endcsname
  \expandafter\let\expandafter\bibitemShut\csname NAT@b@shut@#2\endcsname
  \@ifnum{\NAT@merge>\@ne}{%
   \NAT@bibitem@first@sw{%
    \@firstoftwo
   }{%
    \@ifundefined{NAT@b*@#2}{%
     \@firstoftwo
    }{%
     \expandafter\def\expandafter\NAT@num\expandafter{\the\c@NAT@ctr}%
     \@secondoftwo
    }%
   }%
  }{%
   \@firstoftwo
  }%
  {%
   \global\advance\c@NAT@ctr\@ne
   \@ifx{\NAT@tmp\@empty}{\@firstoftwo}{%
    \@secondoftwo
   }%
   {%
    \expandafter\def\expandafter\NAT@num\expandafter{\the\c@NAT@ctr}%
    \global\NAT@stdbsttrue
   }{}%
   \bibitem@fin
   \item[\href{#3}{\hfil\NAT@anchor{#2}{\NAT@num}}]
   \global\let\NAT@bibitem@first@sw\@secondoftwo
   \NAT@bibitem@init
  }%
  {%
   \NAT@anchor{#2}{}%
   \NAT@bibitem@cont
   \bibitem@fin
  }%
  \@ifx{\NAT@tmp\@empty}{%
    \NAT@wrout{\the\c@NAT@ctr}{}{}{}{#2}%
  }{%
    \expandafter\NAT@ifcmd\NAT@tmp(@)(@)\@nil{#2}%
  }%
}
\makeatother

\author{Tobias Black\thanks{Institut f\"ur Mathematik, Universit\"at Paderborn, Warburger Str. 100, 33098 Paderborn, Germany; email: \mbox{tblack@math.upb.de}}, Chunyan Wu\thanks{School of Mathematical Sciences, University of Electronic Science and Technology of China, 611731 Chengdu, China; email: \mbox{wcypde@163.com} }}
\title{Prescribed signal concentration on the boundary: {W}eak solvability in a chemotaxis-{S}tokes system with proliferation}
\date{}

\begin{document}
\maketitle
\begin{abstract}
\noindent
{\textbf{Abstract:}
We study a chemotaxis-Stokes system with signal consumption and logistic source terms of the form
\noindent 
\begin{align*}
\left\{
\begin{array}{r@{\ }l@{\quad}l@{\quad}l@{\,}c}
n_{t}+u\cdot\!\nabla n&=\Delta n-\nabla\!\cdot(n\nabla c)+\kappa n-\mu  n^{2},\ &x\in\Omega,& t>0,\\
c_{t}+u\cdot\!\nabla c&=\Delta c-nc,\ &x\in\Omega,& t>0,\\
u_{t}&=\Delta u+\nabla P+n\nabla\phi,\ &x\in\Omega,& t>0,\\
\nabla\cdot u&=0,\ &x\in\Omega,& t>0,\\
\big(\nabla n-n\nabla c\big)\cdot\nu&=0,\quad c=c_{\star}(x),\quad u=0, &x\in\partial\Omega,& t>0,
\end{array}\right.
\end{align*}
where $\kappa\geq0$, $\mu>0$ and, in contrast to the commonly investigated variants of chemotaxis-fluid systems, the signal concentration on the boundary of the domain $\Omega\subset\mathbb{R}^N$ with $N\in\{2,3\}$, is a prescribed time-independent nonnegative function $c_{\star}\in C^{2}\!\big(\overline{\Omega}\big)$.\\
Making use of the boundedness information entailed by the quadratic decay term of the first equation, we will show that the system above has at least one global weak solution for any suitably regular triplet of initial data.
}\\[0.08cm]


\end{abstract}


\newpage
\section{Introduction}\label{sec1:intro}
\noindent Chemotaxis, the oriented movement of bacteria and cells in response to a chemical substance in their surrounding environment, is an important motility scheme in nature. An interesting facet of colonies of such chemotactically active bacteria and cells consists of the possibility to spontaneously generate spatial patterns, as not only witnessed by the experimental findings on the aerobic \emph{Bacillus subtilis} (\cite{fujikawaFractalGrowthBacillus1989,matsushitaDiffusionlimitedGrowthBacterial1990,Dombrowski04}) but also in settings where the attracting signal is produced by the cells themselves (\cite{KS70,woodwardSpatiotemporalPatternsGenerated1995}).
This emergence of spatial structures, captivating biologists and mathematicians alike, has lead to an intensive study of chemotaxis systems in the past decades and is still garnering attention in the field of mathematical modeling and analysis. (See also the surveys \cite{Ho03,BBWT15,LanWin_JDMV_19}.)

In order to study the plume-like aggregation patterns observed to occur when a population of \emph{Bacillus subtilis} is suspended in a drop of water, the authors of \cite{tuval2005bacterial} proposed a model of the form
\begin{align}\label{ctns}
\left\{
\begin{array}{r@{\ }l@{\quad}l@{\quad}l@{\,}c}
n_{t}+u\cdot\!\nabla n&=\Delta n-\nabla\!\cdot(n\nabla c),\ &x\in\Omega,& t>0,\\
c_{t}+u\cdot\!\nabla c&=\Delta c-nc,\ &x\in\Omega,& t>0,\\
u_{t}+(u\cdot\nabla)u&=\Delta u+\nabla P+n\nabla\phi,\ &x\in\Omega,& t>0,\\
\nabla\cdot u&=0,\ &x\in\Omega,& t>0,
\end{array}\right.
\end{align}
where $n,c,u,P$ denote the density of the bacteria, the oxygen concentration, the velocity field of the incompressible fluid and the associated pressure, respectively, $\phi$ is a prescribed gravitational potential and $\Omega$ is a bounded domain in $\R^\dimN$. While the authors of \cite{tuval2005bacterial} suggest to augment the system with a non-zero Dirichlet boundary condition for the chemical at the stress-free fluid-air interface and a no-flux condition for the bacteria (in fact they even propose mixed boundary conditions distinguishing between the bottom layer of the drop and the fluid-air interface), a large part of the literature on chemotaxis-fluid systems only considers no-flux conditions for both $n$ and $c$ and a no-slip condition for $u$.

In this setting the global solvability of \eqref{ctns} is well studied and most of the remaining problems remain in the case of $\dimN=3$. Actually, for $\dimN=2$ global classical solutions and their stabilization properties have been established in \cite{win_fluid_CPDE12} and \cite{win-stab2d-ArchRatMechAna12}, respectively. Whereas, in the higher dimensional setting it was shown in \cite{win_globweak3d-AHPN16,win_chemonavstokesfinal_TransAm17} that \eqref{ctns} possesses at least one global weak solution, which becomes smooth after some possibly large waiting time. A recent study by the same author also reveals that on small time-scales (possible) singularities can only arise in a set of measure zero (\cite{win20-leray-structure-ct-fluid}). Similar results have also been established in models where the bacteria are assumed to obey a logistic population growth (i.e. including the term $+\kappa n-\mu n^2$ on the right hand side of the first equation). In fact, existence of weak solutions was shown in \cite{Vorotnikov-WeakSol-CMS14} and \cite{Lan-Longterm_M3AS16} considers the eventual smoothness of weak solutions in 3D. Analytical results providing pattern formation as discovered in the experiments, however, are still missing, which raised the question whether the assumed boundary conditions should be adjusted for further advances. 

Under consideration of different boundary conditions, the knowledge of \eqref{ctns} is quite enigmatic, with most of the current results on existence theory only discussing the two-dimensional setting or relying on the inclusion of small changes to \eqref{ctns}, like logistic growth terms, an enhanced diffusion rate for the bacteria or the consideration of Stokes fluid (i.e. dropping $(u\cdot\nabla)u$ in the third equation) and even then solutions can often only be obtained with quite mild regularity. In this regard, the work \cite{braukhoffGlobalWeakSolution2017} contains the most intricate result in this direction, with the treatment of \eqref{ctns} with logistic growth terms under the Robin boundary condition $\frac{\partial c}{\partial\nu}=1-c$ on $\romega$. The author proves the existence of global classical solutions in 2D and global weak solutions in 3D. Additional results featuring a Robin boundary condition in fluid-free (i.e. $u\equiv0$) variants of \eqref{ctns} have been investigated in \cite{braukhoffStationarySolutionsChemotaxisconsumption2019} and \cite{fuestLongtermBehaviourParabolic2021}. The former considers a stationary (and hence doubly elliptic) system and establishes existence and uniqueness of a classical solution for any prescribed mass $M:=\intomega n>0$. The latter studies a parabolic-elliptic variant and attains results on  global and bounded classical solutions and their long-term behavior. The recent result in \cite{wuAsymptoticDynamicsChemotaxisNavier2020} provides the existence of global weak solutions to the two dimensional version of \eqref{ctns} with superlinear diffusion (i.e. replacing $\Delta n$ by $\Delta n^m$ with $m>1$ in the first equation) and Robin boundary condition for $c$. Concerning non-zero Dirichlet data for $c$ we are only aware of two unpublished works. The first proves global existing generalized solutions in 3D for the Stokes variant of \eqref{ctns} (\cite{WWX2019-local_energy_estimates_prescribed_signal}) and the second provides global generalized solutions for $\dimN\geq2$ in a Stokes variant of \eqref{ctns} with nonlinear diffusion satisfying $m\geq 1$ for $\dimN=2$ and $m>\frac{3\dimN-2}{2\dimN}$ if $\dimN\geq3$ (\cite{WWX2019-global_mass-preserving_dirichlet_signal}). Results on more regular solutions and included logistic population growth appear to be missing for the Dirichlet boundary data case.

(See also \cite{lorz-M3AS10} and \cite{pengGlobalSolutionsCoupled2018,pengGlobalExistenceConvergence2019} for first analytical results concerning well-posedness of systems closely related to \eqref{ctns} with mixed boundary conditions, \cite{knosallaGlobalSolutionsAerotaxis2017} for a more general fluid-free one-dimensional system with non-zero Dirichlet or Neumann boundary data and \cite{tuval2005bacterial,chertockSinkingMergingStationary2012,leeNumericalInvestigationFalling2015} for numerical studies related to \eqref{ctns}.)

\noindent\textbf{Main results.} 
Motivated by the observations above, we are going to consider a chemotaxis-Stokes system with logistic population growth of the form
\begin{align}\label{ct-log-dir}
\left\{
\begin{array}{r@{\ }l@{\quad}l@{\quad}l@{\,}c}
n_{t}+u\cdot\!\nabla n&=\Delta n-\nabla\!\cdot(n\nabla c)+\kappa n-\mu  n^{2},\ &x\in\Omega,& t>0,\\
c_{t}+u\cdot\!\nabla c&=\Delta c-nc,\ &x\in\Omega,& t>0,\\
u_{t}&=\Delta u+\nabla P+n\nabla\phi,\ &x\in\Omega,& t>0,\\
\nabla\cdot u&=0,\ &x\in\Omega,& t>0,\\
\big(\nabla n-n\nabla c\big)\cdot\nu&=0,\quad c=\cs(x),\quad u=0, &x\in\romega,& t>0,\\
n(\cdot,0)&=n_0,\quad c(\cdot,0)=c_0,\quad u(\cdot,0)=u_0 &x\in\Omega,&
\end{array}\right.
\end{align}
in a smoothly bounded domain $\Omega\subset\R^\dimN$ with $\dimN\in\{2,3\}$ and $\nu$ denoting the outward normal vector field on $\romega$. We prescribe $\kappa\geq0$, $\mu>0$, a time constant function $\cs$ satisfying
\begin{align}\label{eq:cstar-reg}
&\cs\in\CSp{2}{\bomega}\quad\text{with}\quad \cs\geq0,
\end{align}
a gravitational potential function $\phi$ fulfilling
\begin{align}\label{phi-def}
\phi\in\W[2,\infty]
\end{align}
and initial data $(n_0,c_0,u_0)$ satisfying
\begin{align}\label{IR} 
\left\{\begin{array}{r@{\,}c@{\,}l@{\ \,}l}
 n_0&\in&\CSp{0}{\bomega}&\text{is nonnegative with }\ n_0\not\equiv0,\\
 c_0&\in&\W[1,q]&\text{is positive in }\Omega\text{ with}\ c_{0}=\cs\text{ on }\romega,\\
 u_0&\in&\DA&
\end{array}\right. 
\end{align}
with $q>\dimN$, $\fracpow\in(\frac{\dimN}{4},1)$. Herein, $A:=-\mathcal{P} \Delta$ denotes the Stokes operator with its domain $D(A):=W^{2,2}\left(\Omega ; \R^{\dimN}\right) \cap W_{0}^{1,2}\left(\Omega ; \R^{\dimN}\right) \cap L_{\sigma}^{2}(\Omega)$ with $L_{\sigma}^{2}(\Omega):=\left\{\varphi\in\LSp{2}{\Omega;\R^\dimN}\,\vert\,\nabla\cdot\varphi=0\right\}$ and $\mathcal{P}$ stands for the Helmholtz projection of $\LSp{2}{\Omega;\R^{\dimN}}$ onto $L_{\sigma}^{2}(\Omega).$\smallskip

\begin{theorem}\label{theo:1} 
Let $\dimN\in\{2,3\}$ and $\Omega\subset\R^\dimN$ be a bounded domain with smooth boundary. Suppose that $\kappa\geq0$ and $\mu>0$ and that the functions $\cs$ and $\phi$ satisfy \eqref{eq:cstar-reg} and \eqref{phi-def}, respectively. Then, for any $n_0,c_0$ and $u_0$ complying with \eqref{IR}, the system \eqref{ct-log-dir} admits at least one global weak solution $(n,c,u)$ in the sense of Definition \ref{def:sol}.
\end{theorem}

\noindent\textbf{Outline.}
In Section \ref{sec2:weak-sol} we will recall the definition of a global weak solution. Section \ref{sec3:globex} will be devoted to the introduction of families of appropriately regularized systems and their time-global classical solvability. On the path toward time-global classical solvability of the approximating system, we will also establish a first set of basic a priori estimates. The commonly employed testing procedures in chemotaxis systems, however, rely heavily on the Neumann boundary conditions and hence adjustments in the treatment of $c$ are necessary here. The substantial regularity information on $n$, as entailed by the quadratic decay present in the first equation, will be the driving force for the distillation of bounds on the gradient of $c$ (see Lemma \ref{lem:diff-ineq-c}), which are an important cornerstone of our further analysis. In Section \ref{sec4:a-priori-n} we will concern ourselves with improving the bounds on $n$, where, in particular, time-space information on $\nabla n$ is the main objective of the section. In Section \ref{sec5:time-reg} we prepare estimates on the time-derivatives, which upon combination with boundedness results of previous sections allows for the construction of a limit object by means of an Aubin--Lions type argument at the start of Section \ref{sec6:limit}. Finally, in the second part of Section \ref{sec6:limit}, we will verify that the limit solution indeed satisfies the properties required of a global weak solution.

\setcounter{equation}{0} 
\section{Definition of global weak solutions}\label{sec2:weak-sol}
Before we start with our analysis let us briefly recount the necessary properties for a global weak solution in the following definition, where here and below we set $W_{0,\sigma}^{1,1}\big(\Omega;\R^\dimN\big):=W_0^{1,1}\big(\Omega;\R^\dimN\big)\cap L_\sigma^2(\Omega)$. 
\begin{definition}\label{def:sol}
A triple $(n,c,u)$ of functions
\begin{align*}
n&\in\LSploc{2}{\bomega\times[0,\infty)}\cap\LSploc{1}{[0,\infty);\W[1,1]},\\ c&\in\LSploc{1}{[0,\infty);\W[1,1]}\quad\text{with}\quad c-\cs\in\LSplocb{1}{[0,\infty);W_0^{1,1}(\Omega)},\\ u&\in\LSplocb{1}{[0,\infty);W_{0,\sigma}^{1,1}\big(\Omega;\R^\dimN\big)}
\end{align*}
with $n\geq0$ and $c\geq0$ in $\bomega\times[0,\infty)$, will be called a global weak solution of \eqref{ct-log-dir} if
\begin{align*}
&nc\;\text{ belongs to }\LSploc{1}{\bomega\times[0,\infty)},\\
\text{if}\;\  n\nabla c,\;\ &nu\;\text{ and }\; cu\;\text{ belong to }\LSploc{1}{\bomega\times[0,\infty);\R^\dimN},
\end{align*}
if
\begin{align}\label{eq:def-eq-n}
-\intinfomega& n\varphi_t-\intomega n_0\varphi(\cdot,0)\\=\,&-\intinfomega\nabla n\cdot\nabla\varphi+\intinfomega n (\nabla c\cdot\nabla\varphi)+\kappa\intinfomega n\varphi-\mu\intinfomega n^2\varphi+\intinfomega n (u\cdot\nabla\varphi)\nonumber
\end{align}
holds for all $\varphi\in C_0^\infty\big(\bomega\times[0,\infty)\big)$, if
\begin{align}\label{eq:def-eq-c}
-\intinfomega c\hat\varphi_t-\intomega c_0\hat\varphi(\cdot,0)=-\intinfomega\nabla c\cdot\nabla\hat\varphi-\intinfomega nc\hat\varphi+\intinfomega c (u\cdot\nabla\hat\varphi)
\end{align}
is valid for all $\hat{\varphi}\in C_0^\infty\big(\Omega\times[0,\infty)\big)$, and if
\begin{align}\label{eq:def-eq-u}
-\intinfomega u\cdot\psi_t-\intomega u_0\cdot\psi(\cdot,0)=-\intinfomega\nabla u\cdot\nabla\psi+\intinfomega n(\nabla\phi\cdot\psi)
\end{align}
is fulfilled for all $\psi\in C_{0}^\infty\big(\Omega\times[0,\infty);\R^\dimN\big)$ with $\nabla\cdot\psi\equiv0$.
\end{definition}

\setcounter{equation}{0} 
\section{Global existence of approximate solutions and essential regularity estimates}\label{sec3:globex}
The global weak solution asserted by Theorem \ref{theo:1} will be obtained as a limit object of solutions to certain regularized problems. To this end, for a fixed family $(\rhoep)_{\epsi\in(0,1)}\subset C_0^\infty(\Omega)$ of smooth cut-off functions satisfying
\begin{align*}
0\leq\rhoep(x)\leq 1\quad\text{for all }x\in\Omega\quad\text{such that }\rhoep\nearrow1\text{ as }\epsi\searrow0,
\end{align*}
we introduce the corresponding family of approximating problems to \eqref{ct-log-dir} given by
\begin{align}\label{approxprob}
\left\{
\begin{array}{r@{\ }l@{\quad}l@{\quad}l@{\,}c}
n_{\epsi t}+\uep\cdot\!\nabla \nep&=\Delta \nep-\nabla\!\cdot\big(\rhoep\fep(\nep)\nep\nabla \cep\big)+\kappa \nep-\mu  \nep^{2},\ &x\in\Omega,& t>0,\\
c_{\epsi t}+\uep\cdot\!\nabla \cep&=\Delta \cep-\gep(\nep)\cep,\ &x\in\Omega,& t>0,\\
u_{\epsi t}&=\Delta \uep+\nabla P_\epsi+\nep\nabla\phi,\quad\nabla\cdot\uep=0,\ &x\in\Omega,& t>0,\\
\frac{\partial\nep}{\partial\nu}&=0,\qquad\ \cep=\cs(x),\quad \qquad \uep=0, &x\in\romega,& t>0,\\
\nep(\cdot,0)&=n_0,\quad \cep(\cdot,0)=c_0,\quad \uep(\cdot,0)=u_0 &x\in\Omega,&
\end{array}\right.
\end{align}
where $\fep(s):=\frac{1}{(1+\epsi s)^3}$ and $\gep(s):=\frac{s}{1+\epsi s}$ for $s\geq0$ and $\epsi\in(0,1)$.

Due to the non-homogeneous boundary condition, this form of the second equation of \eqref{approxprob}, however, is not easily accessible for Dirichlet heat semigroup estimates we will draw on in our following analysis and hence, we substitute $\hcep:=\cs-\cep$ and accordingly rewrite the system into the equivalent formulation
\begin{align}\label{transf-approxprob}
\left\{
\begin{array}{r@{\ }l@{\quad}l@{\quad}l@{\,}c}
n_{\epsi t}+\uep\cdot\!\nabla \nep&=\Delta \nep+\nabla\!\cdot\big(\rhoep\fep(\nep)\nep(\nabla \hcep-\nabla \cs)\big)+\kappa \nep-\mu \nep^2,\ &x\in\Omega,& t>0,\\
{\hat c}_{\epsi t}+\uep\cdot\!\nabla \hcep&=\Delta \hcep-\gep(\nep)\hcep-\Delta\cs+\gep(\nep)\cs+\uep\cdot\nabla\cs,\ &x\in\Omega,& t>0,\\
u_{\epsi t}&=\Delta \uep+\nabla P_\epsi+\nep\nabla\phi,\quad\nabla\cdot\uep=0,\ &x\in\Omega,& t>0,\\
\frac{\partial\nep}{\partial\nu}&=0,\qquad\ \hcep=0,\quad \qquad \uep=0, &x\in\romega,& t>0,\\
\nep(\cdot,0)&=n_0,\quad \hcep(\cdot,0)=\cs-c_0,\quad \uep(\cdot,0)=u_0 &x\in\Omega,&
\end{array}\right.
\end{align}
where, in light of the assumed regularity of $\cs$, all important properties can be easily transferred back to \eqref{approxprob}. The transformed system will only play a role in the proof of time local existence of solutions (Lemma~\ref{lem:locex}) and in the proof that the maximal existence time for fixed $\epsi\in(0,1)$ is actually infinite (Lemma~\ref{lem:globex}), as otherwise our analysis in the latter will not necessarily require semigroup arguments for the second component of the systems.

Now, let us begin by establishing time-local existence of solutions to \eqref{transf-approxprob} (and in turn \eqref{approxprob}) by means of well-established fixed point arguments.

\begin{lemma}
\label{lem:locex}
Let $\Omega\subset \R^\dimN$ be a bounded domain with smooth boundary, $q>\dimN$, $\fracpow\in(\frac{\dimN}{4},1)$, $\kappa\geq0$, $\mu>0$. Suppose that $\cs$ and $\phi$ satisfy \eqref{eq:cstar-reg} and \eqref{phi-def}, respectively, and that $n_0$, $c_0$ and $u_0$ comply with \eqref{IR}. Then for any $\epsi\in(0,1)$, there exist $\Tme\in(0,\infty]$ and a uniquely determined triple $(\nep,\cep,\uep)$ of functions 
\begin{align*}
\nep&\in\CSp{0}{\bomega\times[0,\Tme)}\cap\CSp{2,1}{\bomega\times(0,\Tme)},\\
\cep&\in\CSp{0}{\bomega\times[0,\Tme)}\cap\CSp{2,1}{\bomega\times(0,\Tme)}\cap\LSplocb{\infty}{(0,\Tme);\W[1,q]},\\
\uep&\in\CSp{0}{\bomega\times[0,\Tme);\R^\dimN}\cap\CSp{2,1}{\bomega\times(0,\Tme);\R^\dimN},
\end{align*}
which, together with some $P_\epsi\in\CSp{1,0}{\bomega\times(0,\Tme)}$, solve \eqref{approxprob} in the classical sense and satisfy $\nep\geq0$ and $\cep\geq0$ in $\bomega\times[0,\Tme)$. Moreover, either $\Tme=\infty$ or
\begin{align}\label{eq:extensibility}
\limsup_{t\nearrow\Tme}\left(\|\nep(\cdot,t)\|_{\Lo[\infty]}+\|\cep(\cdot,t)\|_{\W[1,q]}+\|A^\fracpow\uep(\cdot,t)\|_{\Lo[2]}\right)=\infty.
\end{align}
\end{lemma}

\begin{bew}
Augmenting well-established fixed point arguments as e.g. presented in \cite[Lemma 2.1]{Win-ct_fluid_3d-CPDE15} and \cite[Lemma 3.1]{BBWT15} we will first establish time-local existence for the transformed system \eqref{transf-approxprob}, which afterwards, in view of the substitution $\cep=\cs-\hcep$, can be easily transferred back to the corresponding statement for \eqref{approxprob}. For the sake of completeness let us specify the main steps involved:

First, for some large $R>0$ and $T\in(0,1]$, to be specified later, we define the Banach space $X:=\LSpb{\infty}{(0, T) ; \CSp{0}{\bomega} \times W_0^{1,q}(\Omega) \times \DA}$ and its subset 
\begin{align*}
S:=\Big\{(\nep,\hcep,\uep)\in X\,\big\vert\,\|\nep(\cdot,t)\|_{\Lo[\infty]}+\|\hcep(\cdot,t)\|_{\W[1,q]}+\|A^\fracpow \uep(\cdot,t)\|_{\Lo[2]}\leq R\ \text{for a.e. }t\in(0,T)\Big\}.
\end{align*}
Next, denoting by $\big(e^{t\Delta}\big)_{t\geq0}$, $\big(e^{t\Delta'}\big)_{t\geq0}$ and $\big(e^{-tA}\big)_{t\geq0}$ the Neumann heat semigroup, the Dirichlet heat semigroup and the Stokes semigroup with Dirichlet boundary data, respectively, we utilize introduce the mapping $\Phi:=(\Phi_1,\Phi_2,\Phi_3):X\to X$ given by
\begin{align}\label{eq:locex-n}
&\Phi_1(\nep,\hcep,\uep)(\cdot,t)\nonumber\\
:=\ &e^{t\Delta}n_0+\intot e^{(t-s)\Delta}\Big(\nabla\cdot\big(-\uep\nep+\rhoep\fep(\nep)\nep(\nabla \hcep-\nabla \cs)\big)+\kappa \nep-\mu \nep^2\Big)(\cdot,s)\intd s,\\
&\Phi_2(\nep,\hcep,\uep)(\cdot,t)\nonumber\\
:=\ &e^{t\Delta'}(\cs-c_0)+\intot e^{(t-s)\Delta'}\big(\uep\cdot\nabla(\cs-\hcep)+\gep(\nep)(\cs-\hcep)-\Delta\cs\big)(\cdot,s)\intd s,\label{eq:locex-c}
\intertext{and}
&\Phi_3(\nep,\hcep,\uep)(\cdot,t):=e^{-tA}u_0+\intot e^{-(t-s)A}\HP\big(\nep\nabla\phi\big)(\cdot,s)\intd s\quad\text{for }t\in(0,T).\label{eq:locex-u}
\end{align}
We will now show that $\Phi$ acts as a contracting self map on $S$, provided $R$ and $T$ are suitably fixed beforehand. 
Dropping the $\epsi$-subscript for readability, we pick $(n_1,\hc_1,u_1)$, $(n_2,\hc_2,u_2)\in S$ and observe that according to \eqref{eq:locex-n}
\begin{align*}
&\Big\|\big(\Phi_1(n_1,\hc_1,u_1)-\Phi_1(n_2,\hc_2,u_2)\big)(\cdot,t)\Big\|_{\Lo[\infty]}\\
\leq\ &\intot\Big\|e^{(t-s)\Delta}\nabla\cdot\big(-n_1(u_1-u_2)-u_2(n_1-n_2)\big)(\cdot,s)\Big\|_{\Lo[\infty]}\intd s\\
&+\intot\Big\|e^{(t-s)\Delta}\nabla\cdot\big(\rhoep\fep(n_1)(\nabla\hc_1-\nabla\cs)(n_1-n_2)\big)(\cdot,s)\Big\|_{\Lo[\infty]}\intd s\\
&+\intot\bigg\|e^{(t-s)\Delta}\nabla\cdot\Big(\rhoep \fep(n_1)n_2(\nabla\hc_1-\nabla\hc_2)+\rhoep n_2(\nabla\hc_2-\nabla\cs)\big(\fep(n_1)-\fep(n_2)\big)\Big)(\cdot,s)\bigg\|_{\Lo[\infty]}\intd s\\
&+\intot\Big\|e^{(t-s)\Delta}\big(\kappa(n_1-n_2)-\mu(n_1+n_2)(n_1-n_2)\big)(\cdot,s)\Big\|_{\Lo[\infty]}\intd s\quad\text{for }t\in(0,T).
\end{align*}
Hence, drawing on semigroup estimates as e.g. provided by \cite[Lemma 1.3]{win10jde}, \cite[Lemma 2.1]{caolan16_smalldatasol3dnavstokes} and \cite[Lemma 3.1]{Lan15-threshold} we can find $C_1=C_1(\Omega)>0$ such that
\begin{align*}
&\Big\|\big(\Phi_1(n_1,\hc_1,u_1)-\Phi_1(n_2,\hc_2,u_2)\big)(\cdot,t)\Big\|_{\Lo[\infty]}\\
\leq\ &C_1\intot\big(1+(t-s)^{-\frac{1}{2}}\big)\big(\|n_1\|_{\Lo[\infty]}\big\|u_1-u_2\|_{\Lo[\infty]}+\|u_2\|_{\Lo[\infty]}\|n_1-n_2\|_{\Lo[\infty]}\big)(s)\intd s\\
&+C_1\intot\big(1+(t-s)^{-\frac{1}{2}-\frac{\dimN}{2q}}\big)\big(\big\|\nabla\hc_1-\nabla\cs\big\|_{\Lo[q]}\big\|n_1-n_2\big\|_{\Lo[\infty]}\big)(s)\intd s\\
&+C_1\intot\big(1+(t-s)^{-\frac{1}{2}-\frac{\dimN}{2q}}\big)\big(\|n_2\|_{\Lo[\infty]}\|\nabla\hc_1-\nabla\hc_2\|_{\Lo[q]}\big)(s)\intd s\\
&+C_1\intot\big(1+(t-s)^{-\frac{1}{2}-\frac{\dimN}{2q}}\big)\big(\|n_2\|_{\Lo[\infty]}\|\nabla\hc_2-\nabla\cs\|_{\Lo[q]}\|\fep(n_1)-\fep(n_2)\|_{\Lo[\infty]}\big)(s)\intd s\\
&+\kappa\intot\big\|n_1(\cdot,s)-n_2(\cdot,s)\big\|_{\Lo[\infty]}\intd s+\mu\intot\big(\|n_1+n_2\|_{\Lo[\infty]}\|n_1-n_2\|_{\Lo[\infty]}\big)(s)\intd s\quad\text{for }t\in(0,T),
\end{align*}
where we also used the facts that $\rhoep\leq 1$ in $\Omega$, $\fep\leq 1$ in $[0,\infty)$. Moreover, we have $|\fep(a)-\fep(b)|\leq |a-b||a^2+b^2+ab+3a+3b+3|$ for $a,b\in[0,\infty)$ and all $\epsi\in(0,1)$ and $q>\dimN$ as well as $D(A^\fracpow)\hookrightarrow\CSp{\theta}{\bomega}$ for any $\theta\in(0,2\fracpow-\tfrac{\dimN}{2})$ (e.g. \cite[Lemma~III.2.4.3]{sohr} and \cite[Thm.~5.6.5]{evans}) so that we can find $C_2=C_2(\cs,\kappa,\mu,\fracpow,\dimN,q,R,\Omega)$ such that
\begin{align}\label{eq:loc-ex-est-phi1}
\sup_{t\in(0,T)}\Big\|\big(\Phi_1(n_1,\hc_1,u_1)&-\Phi_1(n_2,\hc_2,u_2)\big)(\cdot,t)\Big\|_{\Lo[\infty]}\nonumber\\
&\leq C_2 \big(T+T^{\frac{1}{2}}+T^{\frac{1}{2}-\frac{\dimN}{2q}}\big)
\big\|(n_1-n_2,\hc_1-\hc_2,u_1-u_2)\big\|_{X}.
\end{align}
Similarly, noting that $|\gep(a)-\gep(b)|\leq |a-b|$ for $a,b\in[0,\infty)$ we can draw on semigroup theory for the Dirichlet heat semigroup (see \cite[Proposition~48.4]{QS07} and \cite{hen81}) and \eqref{eq:locex-c} to conclude the existence of $C_3=C_3(\cs,\fracpow,\dimN,q,R,\Omega)>0$ satisfying
\begin{align}\label{eq:loc-ex-est-phi2}
&\sup_{t\in(0,T)}\Big\|\big(\Phi_2(n_1,\hc_1,u_1)-\Phi_2(n_2,\hc_2,u_2)\big)(\cdot,t)\Big\|_{\W[1,q]}\nonumber\\
\leq\ &\sup_{t\in(0,T)}\intot \Big\|e^{(t-s)\Delta'}\big((u_1-u_2)\cdot\nabla(\cs-\hc_1)-u_2\cdot\nabla(\hc_1-\hc_2)\big)(\cdot,s)\Big\|_{\W[1,q]}\intd s\nonumber\\
&\hspace*{1.8cm}+\sup_{t\in(0,T)}\intot \Big\|e^{(t-s)\Delta'}\Big(\big(\gep(n_1)-\gep(n_2)\big)(\cs-\hc_1)-\gep(n_2)(\hc_1-\hc_2)\Big)(\cdot,s)\Big\|_{\W[1,q]}\intd s\nonumber\\
\leq\ &C_3\big(T+T^\frac{1}{2}\big)\big\|(n_1-n_2,\hc_1-\hc_2,u_1-u_2)\big\|_X.
\end{align}
For \eqref{eq:locex-u} we rely on semigroup estimates for the Stokes equation (cf. \cite[Lemma 2.3]{caolan16_smalldatasol3dnavstokes} and \cite[Lemma 3.1]{Win-ct_fluid_3d-CPDE15}) to obtain $C_4=C_4(\phi,\fracpow,\dimN,R,\Omega)>0$ such that
\begin{align}\label{eq:loc-ex-est-phi3}
\sup_{t\in(0,T)}\Big\|A^\fracpow\big(\Phi_3(n_1,\hc_1,u_1)-\Phi_3(n_2,\hc_2,u_2)\big)(\cdot,t)\Big\|_{\Lo[2]}\leq C_4T^{1-\fracpow}\big\|(n_1-n_2,\hc_1-\hc_2,u_1-u_2)\big\|_{X},
\end{align}
so that collecting \eqref{eq:loc-ex-est-phi1}, \eqref{eq:loc-ex-est-phi2} and \eqref{eq:loc-ex-est-phi3} yields
\begin{align}\label{eq:loc-ex-contraction}
\big\|\Phi(n_1,\hc_1,u_1)-\Phi(n_2,\hc_2,u_2)\big\|_X\!\leq\! C_5\big(T+T^\frac{1}{2}+T^{\frac{1}{2}-\frac{\dimN}{2q}}+T^{1-\fracpow}\big)\big\|(n_1-n_2,\hc_1-\hc_2,u_1-u_2)\big\|_X,
\end{align}
with some $C_5=C_5(\cs,\kappa,\mu,\phi,\fracpow,\dimN,q,R,\Omega)>0$. Moreover, since the Dirichlet heat-semigroup estimates provide $C_6=C_6(\Omega)>0$ such that 
\begin{align*}
\intoT \big\|e^{(t-s)\Delta'}\Delta\cs\big\|_{\W[1,q]}\intd s\leq C_6\big(T+T^\frac12\big)\|\Delta\cs\|_{\Lo[q]},
\end{align*}
we find that for some $C_7=C_7(\cs,\Omega)>0$
\begin{align}\label{eq:loc-ex-selfmap}
\big\|\Phi(n,\hc,u)\big\|_{X}&\leq \big\|\Phi(n,\hc,u)-\Phi(0,0,0)\big\|_X+\big\|\Phi(0,0,0)\big\|_X\nonumber\\
&\leq \big\|\Phi(n,\hc,u)-\Phi(0,0,0)\big\|_X+\big\|(n_0,\cs-c_0,u_0)\big\|_X+C_6\big(T+T^\frac12\big)\|\Delta\cs\|_{\Lo[q]}\\
&\leq C_5\big(T+T^\frac{1}{2}+T^{\frac{1}{2}-\frac{\dimN}{2q}}+T^{1-\fracpow}\big)\big\|(n,\hc,u)\big\|_X+\big\|(n_0,\cs-c_0,u_0)\big\|_X+C_7\big(T+T^\frac12\big).\nonumber
\end{align}
Hence, by first taking $R>3\max\big\{\big\|(n_0,\cs-c_0,u_0)\big\|_{X},2C_7\big\}$ and then $T\in(0,1]$ sufficiently small such that $C_5\big(T+T^\frac{1}{2}+T^{\frac{1}{2}-\frac{\dimN}{2q}}+T^{1-\fracpow}\big)<\frac{1}{3}$, we see from \eqref{eq:loc-ex-selfmap} and \eqref{eq:loc-ex-contraction} that indeed $\Phi$ is a contraction map on $S$ and aided by Banach's fixed point theorem we obtain a unique $(\nep,\hcep,\uep)\in S$ with $\Phi(\nep,\hcep,\uep)=(\nep,\hcep,\uep)$. In light of standard bootstrapping procedures drawing on regularity theories for parabolic equations and the Stokes semigroup \cite{PorzVesp93,Solonnikov2007,LSU} one can verify that $(\nep,\hcep,\uep)$ actually satisfies the claimed regularity properties, which then entails the existence of a corresponding $P_{\epsi}$ such that $(\nep,\hcep,\uep,P_{\epsi})$ solves \eqref{transf-approxprob} classically in $\Omega\times(0,T)$. Uniqueness of $(\nep,\hcep,\uep)$ can be verified by standard $L^2$ testing procedures for the differences of two assumed solutions. Noticing that the choice of $T$ only depends on fixed system parameters and the initial data, we may iterate the arguments (with different initial data and possibly larger $R$) to extend the solution on a maximal time interval $(0,\Tme)$ such that either $\Tme=\infty$ or 
\begin{align*}
\limsup_{t\nearrow \Tme}\big(\|\nep(\cdot,t)\|_{\Lo[\infty]}+\|\hcep(\cdot,t)\|_{\W[1,q]}+\|A^\fracpow \uep(\cdot,t)\|_{\Lo[2]}\big)=\infty.
\end{align*}
Clearly, by substituting $\cep=\cs-\hcep$ (and recalling \eqref{eq:cstar-reg}), we immediately obtain the desired results for \eqref{approxprob}, where, finally, the nonnegativity of $\nep$ and $\cep$ is entailed by two applications of the maximum principle to the first and second equation of \eqref{approxprob}.
\end{bew}
For the remainder of the work we will now assume that $\dimN\in\{2,3\}$, $\Omega\subset\R^\dimN$, $\kappa\geq0$, $\mu>0$, $q>N$, $\fracpow\in(\frac{\dimN}{4},1)$, $\cs,\phi$ satisfying \eqref{eq:cstar-reg} and \eqref{phi-def}, respectively, and initial data $n_0,c_0,u_0$ obeying \eqref{IR} are fixed and, accordingly, for $\epsi\in(0,1)$ denote by $(\nep,\cep,\uep)$ the triple of functions provided by Lemma \ref{lem:locex} and by $\Tme$ the corresponding maximal existence time.

Time-local existence at hand, we can now proceed with a first set of a priori properties obtained by straightforward integration and an application of the maximum principle.
\begin{lemma}
\label{lem:simple-bounds}
There is $C>0$ such that for any $\epsi\in(0,1)$  the solution $(\nep,\cep,\uep)$ of \eqref{approxprob} satisfies
\begin{align*}
\intomega \nep(\cdot,t)\leq C,\qquad \inttminusomega \nep^{2}\leq C\qquad\text{and}\qquad\|\cep(\cdot,t)\|_{\Lo[\infty]}\leq C\quad\text{for all}\quad t\in(0,\Tme).
\end{align*}
\end{lemma}

\begin{bew}
Making use of the fact that $\uep$ is divergence free, by integrating the first equation of \eqref{approxprob} over $\Omega$ and utilizing integration by parts as well as Young's inequality we deduce that for all $\epsi\in(0,1)$
\begin{align*}
&\frac{\intd}{\intd t}\intomega \nep+ \mu\intomega \nep^2=\kappa\intomega \nep\leq\frac{\mu}{2}\intomega \nep^2+\frac{\kappa^2}{2\mu}{|\Omega|}\quad\text{on }(0,\Tme).
\end{align*}
Employing Young's inequality once more to estimate the quadratic term on the left from below we obtain
\begin{align}\label{eq:n1}
\frac{\intd}{\intd t}\intomega \nep+\mu\intomega \nep+\frac{\mu}{4}\intomega \nep^2\leq\frac{\kappa^2}{2\mu}{|\Omega|}
+\mu|\Omega|\quad\text{on }(0,\Tme)\text{ for all }\epsi\in(0,1), 
\end{align}
which, when combined with the nonnegativity of $\nep$ and an ODE comparison argument, implies 
\[\intomega \nep(\cdot,t)\leq C_1:=\max\left\{\intomega n_0,\,\left(\frac{\kappa^2}{2\mu^2}+ 1\right)|\Omega|\right\}\quad\text{for all } t\in(0,\Tme)\text{ and all }\epsi\in(0,1).\]
Furthermore, integration of \eqref{eq:n1} over $\big((t-1)_+,t\big)$ now provides,
\[
\frac{\mu}{4}\inttminusomega\nep^2\leq \intomega\nep\big(\cdot,(t-1)_+\big)+\frac{\kappa^2}{2\mu}{|\Omega|}
+\mu{|\Omega|}\leq C_1+\frac{\kappa^2}{2\mu}{|\Omega|}
+\mu{|\Omega|}.
\]
for all $t\in(0,\Tme)$ and all $\epsi\in(0,1)$. Finally,  by the maximum principle, we instantly obtain that
\[
 \|\cep(\cdot,t)\|_{\Lo[\infty]}\leq \max\left\{\|\cs\|_{\LSp{\infty}{\romega}},\|c_0\|_{\Lo[\infty]}\right\} \quad\text{for all } t\in(0,\Tme)\text{ and }\epsi\in(0,1),
\]
which completes the proof. 
\end{bew}

In order to distill further uniform bounds from the somewhat sparse (yet sufficiently powerful) space-time information on $\nep^2$ provided by Lemma \ref{lem:simple-bounds}, we state the following comparison result for ordinary differential equations. This lemma is copied from \cite[Lemma 3.4]{LL18-ClassSolLogCTSingSens}, whereto we refer the reader for details of the proof.
\begin{lemma}
\label{lem:ode-lemma}
For some $T\in(0,\infty]$ let $y\in\CSp{1}{(0,T)}\cap\CSp{0}{[0,T)}$, $h\in\CSp{0}{[0,T)}$, $h\geq0$, $C>0$, $a>0$ satisfy
\begin{align*}
y'(t)+ay(t)\leq h(t),\qquad \int_{(t-1)_+}^t h(s)\intd s\leq C
\end{align*}
for all $t\in(0,T)$. Then $y\leq y(0)+\frac{C}{1-e^{-a}}$ throughout $(0,T)$.
\end{lemma}

With the comparison lemma above, we can make now turn to obtain some uniform bounds for the third solution component.
\begin{lemma}
\label{lem:u-stokes-bound}
There is $C>0$ such that for any $\epsi\in(0,1)$ the solution $(\nep,\cep,\uep)$ of \eqref{approxprob} satisfies
\begin{align*}
\intomega|\nabla\uep(\cdot,t)|^2\leq C\qquad\text{and}\qquad \intomega|\uep(\cdot,t)|^6\leq C\qquad\text{for all }t\in(0,\Tme).
\end{align*}
\end{lemma}

\begin{bew}
First, we test the third equation in \eqref{approxprob} against $\uep$, integrate by parts over $\Omega$, and employ the Young and Poincar\'e inequalities as well as \eqref{phi-def} to conclude the existence of $C_1>0$ such that for all $\epsi\in(0,1)$
\begin{align}
\frac{1}{2}\frac{\intd}{\intd t}\intomega|\uep|^2+\frac{1}{2}\intomega |\nabla \uep|^2\le C_1 \intomega \nep^2\label{u1}
\end{align} 
is valid on $(0,\Tme)$. Then, again denoting by $\HP$ the Helmholtz projection and by $A$ the Stokes operator, we multiply the projected third equation by $A\uep$ to obtain $C_2>0$ such that for all $\epsi\in(0,1)$ the inequality
\begin{align}
\frac{1}{2}\frac{\intd}{\intd t}\intomega| A^{\frac{1}{2}}\uep|^2+\intomega|A\uep|^2=\intomega\HP[\nep\nabla\phi]\cdot A\uep\leq \frac{1}{2}\intomega|A\uep|^2+C_2\intomega \nep^2 \label{u2}
\end{align}
holds on $(0,\Tme)$, where we once more made use of the boundeness of $\nabla\phi$ and Young's inequality. In light of the Poincar\'e inequality, a combination of \eqref{u1} and \eqref{u2} entails the existence of $C_3,C_4>0$ such that
\begin{align*}
\frac{\intd}{\intd t}\left\{\intomega|\uep|^2+\intomega|\nabla \uep|^2\right\} 
+C_3\left\{\intomega|\uep|^2+\intomega|\nabla \uep|^2\right\}+\intomega|A\uep|^2\leq (2C_1+2C_2)\intomega \nep^2+C_4
\end{align*}
on $(0,\Tme)$ for all $\epsi\in(0,1)$. Drawing on Lemmas  \ref{lem:simple-bounds} and \ref{lem:ode-lemma}, we conclude that there is $C_5>0$ satisfying
\[
\intomega|\nabla\uep(\cdot,t)|^2\leq C_5  \qquad\text{for all }t\in(0,\Tme)\text{ and }\epsi\in(0,1).
\]
Finally, relying on the Sobolev embedding theorem, we find $C_6>0$ such that
\[
\intomega|\uep(\cdot,t)|^6\leq C_6 \qquad\text{for all }t\in(0,\Tme)\text{ and }\epsi\in(0,1).\qedhere
\]
\end{bew}

The uniform bounds on $\uep$ in $\LSp{\infty}{(0,\Tme);\Lo[6]}$ and $\nep$ in $\LSp{2}{\Omega\times(0,\Tme)}$ will now be the key ingredient in obtaining information on $\nabla\cep$. We start by exploiting the fact that $\cs$ is constant in time to establish an ordinary differential inequality for $\intomega|\nabla\cep(\cdot,t)|^2$ on $(0,\Tme)$.

\begin{lemma}
\label{lem:diff-ineq-c}
There exists $C>0$ such that for all $\epsi\in(0,1)$ the solution $(\nep,\cep,\uep)$ of \eqref{approxprob} satisfies
\begin{align*}
\frac{1}{2}\frac{\intd}{\intd t}\intomega|\nabla\cep|^2+\frac{1}{4}\intomega|\Delta\cep|^2\leq C\intomega\nep^2+C
\end{align*}
on $(0,\Tme)$.
\end{lemma}

\begin{bew}
Since the boundary conditions in \eqref{approxprob} imply that $\frac{\partial}{\partial t}c_{\epsi}\big\vert_{\romega}=0$ on $(0,\Tme)$, we can multiply the second equation of \eqref{approxprob} by $-\Delta\cep$ and integrate by parts to find that for all $\epsi\in(0,1)$
\begin{align*}
\frac{1}{2}\frac{\intd}{\intd t}\intomega|\nabla\cep|^2&=-\intomega\Delta\cep c_{\epsi t}+\intromega c_{\epsi t}\frac{\partial\cep}{\partial\nu}=-\intomega|\Delta\cep|^2+\intomega(\uep\cdot\nabla\cep)\Delta\cep+\intomega\gep(\nep)\cep\Delta\cep
\end{align*}
on $(0,\Tme)$. Employing Young's inequality to the last two terms on the right and making use of the fact that $|\gep(s)|\leq s$ for all $s\geq0$ we obtain that for all $\epsi\in(0,1)$
\begin{align}\label{eq:diff-ineq-c-eq1}
\frac{1}{2}\frac{\intd}{\intd t}\intomega|\nabla\cep|^2+\frac{1}{2}\intomega|\Delta\cep|^2\leq\intomega\nep^2\cep^2+\intomega|\uep\cdot\nabla\cep|^2\quad\text{on }(0,\Tme).
\end{align}
In view of Lemma \ref{lem:simple-bounds}, there is $C_1>0$ satisfying
\begin{align}\label{eq:diff-ineq-c-eq2}
\intomega\nep^2\cep^2\leq\|\cep\|_{\Lo[\infty]}^2\intomega\nep^2\leq C_1\intomega\nep^2\quad\text{on }(0,\Tme)\ \text{ for all }\epsi\in(0,1).
\end{align}
Furthermore, relying on the Hölder inequality and Lemma \ref{lem:u-stokes-bound}, we find $C_2>0$ such that for all $\epsi\in(0,1)$ we have
\begin{align*}
\intomega|\uep\cdot\nabla\cep|^2\leq\|\uep\|_{\Lo[6]}^2\|\nabla\cep\|_{\Lo[3]}^2\leq C_2\|\nabla\cep\|_{\Lo[3]}^2\quad\text{on }(0,\Tme).
\end{align*}
The \GNI\ and Lemma \ref{lem:simple-bounds}, moreover, imply the existence of $C_3,C_4>0$ satisfying
\begin{align*}
\|\nabla\cep\|_{\Lo[3]}^2\leq C_3\|\Delta\cep\|_{\Lo[2]}\|\cep\|_{\Lo[6]}+C_3\|\cep\|_{\Lo[6]}^2\leq C_4\|\Delta\cep\|_{\Lo[2]}+C_4
\end{align*}
on $(0,\Tme)$ for all $\epsi\in(0,1)$, so that an application of Young's inequality entails the existence of $C_5>0$ such that for all $\epsi\in(0,1)$ we have
\begin{align}\label{eq:diff-ineq-c-eq3}
\intomega|\uep\cdot\nabla\cep|^2\leq\frac{1}{4}\intomega|\Delta\cep|^2+C_5\quad\text{on }(0,\Tme).
\end{align} 
A combination of \eqref{eq:diff-ineq-c-eq1}-\eqref{eq:diff-ineq-c-eq3} finally shows that with $C:=\max\{C_1,C_5\}$ we obtain
\begin{align*}
\frac{1}{2}\frac{\intd}{\intd t}\intomega|\nabla\cep|^2+\frac{1}{4}\intomega|\Delta\cep|^2\leq C\intomega\nep^2+C
\end{align*}
on $(0,\Tme)$ for all $\epsi\in(0,1)$.
\end{bew}

Next, we combine the recently established differential inequality with the \GNI, the comparison Lemma \ref{lem:ode-lemma} and the space-time bound for $\nep$ from Lemma \ref{lem:simple-bounds} to obtain the following.
\begin{lemma}
\label{lem:grad-c-l4}
There is $C>0$ such that for any $\epsi\in(0,1)$ the solution $(\nep,\cep,\uep)$ of \eqref{approxprob} fulfills
\begin{align*}
\intomega|\nabla\cep(\cdot,t)|^2\leq C,\qquad \inttminusomega|\Delta\cep|^2\leq C\qquad\text{and}\qquad\inttminusomega|\nabla\cep|^4\leq C
\end{align*}
for all $t\in(0,\Tme)$.
\end{lemma}

\begin{bew}
According to the \GNI\ and Lemma \ref{lem:simple-bounds}, there are $C_1,C_2>0$ such that for all $\epsi\in(0,1)$
\[
\|\nabla\cep\|_{\Lo[2]}^2\le C_1\|\Delta\cep\|_{\Lo[2]}\|\cep\|_{\Lo[2]}+C_1\|\cep\|_{\Lo[\infty]}^2\le \frac{1}{4}\|\Delta\cep\|_{\Lo[2]}^2+C_2
\]
on $(0,\Tme)$, which upon combination with the differential inequality for $\cep$ in Lemma \ref{lem:diff-ineq-c}, the bounds obtained in Lemma \ref{lem:simple-bounds} and the ODE-comparison of Lemma \ref{lem:ode-lemma} entails the existence of $C_3>0$ such that for all $\epsi\in(0,1)$ we have
\begin{align}\label{eq:grad-c-l4-eq1}
\intomega|\nabla\cep|^2\le C_3\quad\text{on }(0,\Tme).
\end{align}
Then, returning to the differential inequality for $\cep$ from Lemma \ref{lem:diff-ineq-c}, we obtain $C_4>0$ such that for all $\epsi\in(0,1)$ we have $\inttminusomega|\Delta\cep|^2 \le C_4$ on $(0,\Tme-\tau)$ from straightforward integration of said inequality in light of \eqref{eq:grad-c-l4-eq1} and Lemma \ref{lem:simple-bounds}.
Finally, once again in view of the \GNI, we we find $C_5>0$ such that for all $\epsi\in(0,1)$
\begin{align*}
\int_{(t-1)_+}^t\|\nabla\cep\|_{\Lo[4]}^4\leq C_5\int_{(t-1)_+}^t\|\Delta\cep\|_{\Lo[2]}^2\|\cep\|_{\Lo[\infty]}^2+C_5\int_{(t-1)_+}^t\|\cep\|_{\Lo[\infty]}^4\quad\text{on }(0,\Tme),
\end{align*}
completing the proof by drawing on the previous parts of this lemma and Lemma \ref{lem:simple-bounds}.
\end{bew}

The boundedness property of $\nabla\cep$ in $\LSp{\infty}{(0,\Tme);\Lo[2]}$ was the last missing piece of information necessary for proving time-global existence of solution to \eqref{approxprob}. Augmenting the bounds we established in this Section with additional $\epsi$-dependent bounds in the proof below, we will be able to draw on a Moser--Alikakos-type iteration procedure (see \cite[Lemma A.1]{TaoWin-quasilinear_JDE12}) to finally conclude that for fixed $\epsi$ the maximal existence time $\Tme$ provided by Lemma \ref{lem:locex} is indeed not finite.

\begin{lemma}
\label{lem:globex}
For all $\epsi\in(0,1)$ the solution of \eqref{approxprob} is global in time, i.e. $\Tme=\infty$.
\end{lemma}

\begin{bew}
 We fix $\epsi\in(0,1)$ and assume for contradiction that $T:=\Tme<\infty$. Subsequently we will consider estimates for the quantities appearing in the extensibility criterion \eqref{eq:extensibility} and, as our estimation process relies on employing semigroup arguments to the second component, we will once more return to working in the transformed system \eqref{transf-approxprob}. Since an immediate estimation of $\nep$ in $\LSp{\infty}{\Omega\times(0,T)}$ is out of our reach, we will first establish the boundedness of $\nep$ in $\LSp{\infty}{(0,T);\Lo[6]}$. To this regard, we multiply the first equation of \eqref{transf-approxprob} by $\nep^5$ and integrate by parts to find that due to $\uep$ being divergence free
\begin{align*}
\frac{1}{6}\frac {\intd}{\intd t}\intomega\nep^6+5\intomega\nep^{4}|\nabla\nep|^2
=&-5\intomega\rhoep\fep(\nep)\nep^{5}(\nabla\hcep-\nabla \cs)\cdot\nabla\nep+\kappa\intomega\nep^6-\mu\intomega\nep^7\\
&\quad+\intromega\! \nep^5\big(\nabla\nep-\rhoep\fep(\nep)(\nabla\hcep-\nabla\cs)\big)\cdot\nu-\frac16\!\intromega\!\nep^6(\uep\cdot\nu)
\end{align*}
on $(0,T)$. Here, the last two integrals disappear because of the boundary conditions and the fact that $\rhoep=0$ on $\romega$. Moreover, noticing that $|\fep(s)s^3|\leq\frac{1}{\epsi^3}$ for all $s\geq0$ and that $|\rhoep|\leq 1$ on $\Omega$, we make use of two applications of Young's inequality to find $C_1:=C_1(\kappa,\mu)>0$ such that
\begin{align*}
\frac{1}{6}\frac{\intd}{\intd t}\intomega\nep^6+\frac{5}{2}\intomega\nep^4|\nabla\nep|^2+\intomega\nep^6\leq \frac{5}{2\epsi^6}\intomega|\nabla\hcep-\nabla\cs|^2+C_1\quad\text{on }(0,T).
\end{align*} 
Since $|\nabla\hcep-\nabla\cs|^2=|\nabla\cep|^2$, we conclude from Lemma \ref{lem:grad-c-l4} and a straightforward comparison argument that there is $C_2=C_2(\epsi)>0$ satisfying
\begin{align}\label{eq:n-l6-globex}
\intomega \nep^6\leq C_2\quad\text{on }(0,T).
\end{align}
This bound at hand, we pick $\fracpow\in(\frac{\dimN}{4},1)$ as in \eqref{IR} and then rely on smoothing properties of the Stokes semigroup (e.g. \cite[p.201]{gig86} and \cite[Lemma 3.1]{Win-ct_fluid_3d-CPDE15}), \eqref{IR}, \eqref{phi-def} and \eqref{eq:n-l6-globex} to obtain $C_3=C_3(\epsi)>0$ satisfying
\begin{align*}
\|A^{\fracpow}\uep(\cdot,t)\|_{\Lo[2]}&\leq \|A^\fracpow u_0\|_{\Lo[2]}+\int_0^t\big\|A^\fracpow e^{-(t-s)A}\HP(n(\cdot,s)\nabla\phi)\big\|_{\Lo[2]}\intd s\leq C_3+\frac{C_3 T^{1-\fracpow}}{1-\fracpow}
\end{align*}
for all $t\in(0,T)$, which, due to the embedding $D(A^\fracpow)\hookrightarrow\CSp{\theta}{\bomega}$ for any $\theta\in(0,2\fracpow-\tfrac{\dimN}{2})$ (cf. \cite[Lemma~III.2.4.3]{sohr} and \cite[Thm.~5.6.5]{evans}), also entails that there is some $C_4=C_4(T,\epsi)>0$ such that
\begin{align}\label{eq:u-linfty-globex}
\|\uep(\cdot,t)\|_{\Lo[\infty]}\leq C_4\quad\text{for all }t\in(0,T).
\end{align}
Next, drawing on the Dirichlet heat-semigroup representation of $\hcep$ we find that
\begin{align*}
\nabla\hcep(\cdot,t)&=\nabla e^{t\Delta'}(\cs-c_0)+\int_0^t e^{(t-s)\Delta'}\nabla\big(\uep\cdot\nabla(\cs-\hcep)+\gep(\nep)(\cs-\hcep)-\Delta\cs\big)\quad\text{for all }t\in(0,T).
\end{align*}
Picking $q\in(N+2,6)$ and letting $t_0:=\min\{1,\frac{T}{2}\}$, in light of well-known semigroup estimates (\cite{hen81}), we can hence obtain $C_5>0$ satisfying
\begin{align*}
&\|\nabla\hcep(\cdot,t)\|_{\Lo[q]}\\
\leq\ &C_5\Big(1+t_0^{-\frac12-\frac{\dimN}{2}(\frac{1}{\dimN}-\frac1q)}\Big)\big\|\cs-c_0\big\|_{\Lo[\dimN]}+C_5\int_0^t\!\Big(1+(t-s)^{-\frac{1}{2}-\frac{\dimN}{2}(\frac{1}{2}-\frac{1}{q})}\Big)\Big\|\big(\uep\nabla(\cs-\hcep)\big)(\cdot,s)\Big\|_{\Lo[2]}\!\!\intd s\\
&\quad +C_5\int_0^t\!\Big(1+(t-s)^{-\frac{1}{2}}\Big)\Big\|\big(\nep|\cs-\hcep|+|\Delta\cs|\big)(\cdot,s)\Big\|_{\Lo[q]}\!\intd s\quad\text{for all }t\in(t_0,T),
\end{align*}
where we also relied on the estimate $\gep(s)\leq s$ for $s\geq0$. Here, due to $\cs-\hcep=\cep$, we can make use of \eqref{eq:u-linfty-globex} and Lemma \ref{lem:grad-c-l4} for the first integral and \eqref{eq:n-l6-globex} combined with $q<6$, Lemma~\ref{lem:simple-bounds} and \eqref{eq:cstar-reg} for the second integral, to find $C_6=C_6(T,\epsi)>0$ and, since $\frac{2\dimN}{\dimN-2}\geq 6>q$ entails $-\frac{1}{2}-\frac{\dimN}{2}(\frac{1}{2}-\frac{1}{q})>-1$, also $C_7=C_7(T,\epsi)$ satisfying
\begin{align*}
\|\nabla\hcep(\cdot,t)\|_{\Lo[q]}\leq C_6+C_6\int_0^t\!\Big(2+(t-s)^{-\frac{1}{2}}+(t-s)^{-\frac{1}{2}-\frac{\dimN}{2}(\frac{1}{2}-\frac{1}{q})}\Big)\intd s\leq C_7
\end{align*}
for all $t\in(t_0,T)$. With these bounds we can easily check that a Moser--Alikakos-type iteration procedure (see \cite[Lemma A.1]{TaoWin-quasilinear_JDE12}) becomes applicable to \eqref{transf-approxprob} and that hence there is $C_8=C_8(T,\epsi)>0$ such that $\|\nep\|_{\Lo[\infty]}\leq C_8$ on $(0,T)$. Hence, we find
\begin{align*}
&\limsup_{t\nearrow T} \Big(\|\nep(\cdot,t)\|_{\Lo[\infty]}+\|\cep(\cdot,t)\|_{\W[1,q]}+\|A^\fracpow\uep(\cdot,t)\|_{\Lo[2]}\Big)\\
\leq\ &\limsup_{t\nearrow T} \Big(\|\nep(\cdot,t)\|_{\Lo[\infty]}+\|\hcep(\cdot,t)\|_{\W[1,q]}+\|\cs\|_{\W[1,q]}+\|A^\fracpow\uep(\cdot,t)\|_{\Lo[2]}\Big)<\infty,
\end{align*} 
contradicting \eqref{eq:extensibility} and therefore proving $\Tme=\infty$.
\end{bew}
\setcounter{equation}{0} 
\section{Refined a priori information on \texorpdfstring{$\nep$}{n}}\label{sec4:a-priori-n}
While the uniform bounds for $\cep$ and $\uep$ provided by Section \ref{sec3:globex} would already be strong enough for our limit procedure, we still lack sufficiently good uniform bounds for $\nep$. As it turns out, the space-time bound for $\nabla\cep$ of Lemma \ref{lem:grad-c-l4}, however, can be exploited when considering the functional $y_\epsi(t):=\intomega\big(\nep\ln\nep\big)(\cdot,t)$, which has often been a good resource for information in chemotaxis settings (\cite{Duan2010,Lan-Longterm_M3AS16,win_fluid_CPDE12,win_globweak3d-AHPN16}). We start by formulating a corresponding functional inequality.

\begin{lemma}
\label{lem:diff-ineq-n}
There exists $C>0$ such that for all $\epsi\in(0,1)$ the solution $(\nep,\cep,\uep)$ of \eqref{approxprob} satisfies
\begin{align*}
\frac{\intd}{\intd t}\intomega \nep\ln\nep&+\frac{1}{2}\intomega\frac{|\nabla\nep|^2}{\nep}+\frac{\mu}{2}\intomega\nep^2\ln\nep\leq C\intomega \nep^2+C\intomega|\nabla\cep|^4+C
\end{align*}
on $(0,\infty)$.
\end{lemma}

\begin{bew}
In light of the first equation of \eqref{approxprob}, the fact that $\uep$ is divergence free and two integrations by parts we see that
\begin{align*}
\frac{\intd}{\intd t}\intomega\nep\ln\nep&=\intomega\Big(\nabla\cdot(\nabla\nep-\rhoep\fep(\nep)\nep\nabla\cep)-\uep\cdot\nabla\nep+\kappa\nep-\mu\nep^2\Big)(\ln\nep+1)
\\&=-\intomega\frac{|\nabla\nep|^2}{\nep}+\intomega \rhoep\fep(\nep)(\nabla\nep\cdot\nabla\cep)+\intromega(\ln\nep+1)\big(\nabla\nep-\rhoep\fep(\nep)\nep\nabla\cep)\cdot\nu\\&\hspace*{2.3cm}-\intromega\nep\ln\nep(\uep\cdot\nu)+\kappa\intomega\nep\ln\nep+\kappa\intomega\nep-\mu\intomega\nep^2\ln\nep-\mu\intomega\nep^2
\end{align*}
on $(0,\infty)$ for all $\epsi\in(0,1)$. Observing that again both boundary integrals disappear due to the prescribed boundary conditions and the fact that $\rhoep(x)=0$ for $x\in\romega$ and noting that there is some $C_1>0$ satisfying $\kappa s-\mu s^2\leq C_1$ for all $s\geq0$ and such that $(\kappa s-\frac{\mu}{2}s^2)\ln(s)\leq C_1$ for all $s>0$, we may hence estimate
\begin{align*}
\frac{\intd}{\intd t}\intomega\nep\ln\nep+\intomega\frac{|\nabla\nep|^2}{\nep}+\frac{\mu}{2}\intomega\nep^2\ln\nep\leq\intomega\rhoep\fep(\nep)(\nabla\nep\cdot\nabla\cep)+2C_1
\end{align*}
on $(0,\infty)$ for all $\epsi\in(0,1)$. To further estimate the integral on the right, we make use of the fact that $|\rhoep(x)\fep(s)|=\frac{\rhoep(x)}{(1+\epsi s)^3}\leq 1$ for all $s\geq0$, $x\in\Omega$ and $\epsi\in(0,1)$ and two applications of Young's inequality to obtain
\begin{align*}
\frac{\intd}{\intd t}\intomega\nep\ln\nep+\intomega\frac{|\nabla\nep|^2}{\nep}+\frac{\mu}{2}\intomega\nep^2\ln\nep&\leq\frac{1}{2}\intomega\frac{|\nabla\nep|^2}{\nep}+\frac{1}{2}\intomega\nep|\nabla\cep|^2+2C_1\\
&\leq\frac{1}{2}\intomega\frac{|\nabla\nep|^2}{\nep}+\frac{1}{4}\intomega\nep^2+\frac{1}{4}\intomega|\nabla\cep|^4+2C_1
\end{align*}
on $(0,\infty)$ for all $\epsi\in(0,1)$, which concludes the proof upon the choice of $C:=\max\{\frac{1}{4},2C_1\}$.
\end{bew}

Clearly, we can draw on previously established space-time bounds to extract additional space-time information on $\nabla\sqrt{\nep}$ from the previous Lemma, which in a second interpolation step can also be refined to a bound on $\nabla\nep$ in $\LSp{\frac43}{\Omega\times(0,\infty)}$.

\begin{lemma}
\label{lem:eq-n}
For any $T>0$ there is $C(T)>0$ such that for all $\epsi\in(0,1)$ the solution $(\nep,\cep,\uep)$ of \eqref{approxprob} satisfies
\begin{align*}
\intoTomega\frac{|\nabla\nep|^2}{\nep}\leq C(T)\qquad\text{and}\qquad \intoTomega\nep^2\ln\nep\leq C(T).
\end{align*}
\end{lemma}

\begin{bew}
Integration of the differential inequality featured in Lemma \ref{lem:diff-ineq-n} over $(0,T)$ provides $C_1>0$ such that for all $\epsi\in(0,1)$
\begin{align*}
&\frac{1}{2}\intoTomega\frac{|\nabla\nep|^2}{\nep}+\frac{\mu}{2}\intoTomega\nep^2\ln\nep\\
\leq\ & C_1\intoTomega\nep^2+C_1\intoTomega|\nabla\cep|^4-\intomega\nep(\cdot,T)\ln\nep(\cdot,T)+\intomega\nep(\cdot,0)\ln\nep(\cdot,0)+C_1T.
\end{align*}
Recalling the bounds provided by Lemma \ref{lem:simple-bounds}, Lemma \ref{lem:grad-c-l4} and \eqref{IR} as well as the obvious estimate $-\frac{1}{e}\leq s\ln s$ for all $s\geq 0$, the conclusion is immediate.
\end{bew}

\begin{lemma}
\label{lem:grad-n-time}
For all $T>0$ there exists $C(T)>0$ such that for all $\epsi\in(0,1)$ the solution $(\nep,\cep,\uep)$ of \eqref{approxprob} fulfills
\begin{align*}
\intoTomega|\nabla\nep|^\frac{4}{3}\leq C(T).
\end{align*}
\end{lemma}

\begin{bew}
Rewriting the integral under consideration and employing Young's inequality twice, we find that
\begin{align*}
\intoTomega|\nabla\nep|^\frac{4}{3}
&=\intoTomega\frac{|\nabla \nep||\nabla\nep|^{\frac{1}{3}}\sqrt{\nep}}{\sqrt{\nep}}\\
&\leq \frac{1}{2}\intoTomega\frac{|\nabla\nep|^2}{\nep}+\frac{1}{2} \intoTomega|\nabla\nep|^{\frac{2}{3}}\nep\\
&\leq \frac{1}{2}\intoTomega\frac{|\nabla\nep|^2}{\nep}+\frac{1}{4} \intoTomega|\nabla\nep|^{\frac{4}{3}}+\frac{1}{4}\intoTomega\nep^2
\end{align*}
for all $\epsi\in(0,1)$. Reordering and making use of the bounds provided by Lemmas \ref{lem:simple-bounds} and \ref{lem:eq-n} we obtain the asserted bound.
\end{bew}

\setcounter{equation}{0} 
\section{Regularity estimates for the time derivatives}\label{sec5:time-reg}
As final element for an Aubin--Lions type argument we are going to undertake in Section \ref{sec6:limit}, we now prepare uniform bounds for the time derivatives in suitable spaces.
\begin{lemma}
\label{lem:time-reg}
For every $T>0$ there exists $C(T)>0$ such that for any $\epsi\in(0,1)$ the solution $(\nep,\cep,\uep)$ of \eqref{approxprob} fulfills
\begin{align}\label{eq:n-time-regu}
\int_0^T\|\partial_t\nep\|_{(W_0^{1,4}(\Omega))^*}\leq C(T),
\end{align}
\begin{align}\label{eq:c-time-regu}
\int_0^T\|\partial_t\cep\|_{(W_0^{1,4}(\Omega))^*}\leq C(T),
\end{align}
and 
\begin{align}\label{eq:u-time-regu}
\int_0^T\|\partial_t \uep\|_{(W_{0,\sigma}^{1,2}(\Omega))^*}^{2}\leq C(T).
\end{align}
\end{lemma}

\begin{bew}
Given $T>0$ we fix $\varphi\in\LSpb{\infty}{(0,T);W_0^{1,4}(\Omega)}$ with $\|\varphi\|_{\LSpn{\infty}{(0,T);W_0^{1,4}(\Omega)}}\leq 1$ and test the first equation of \eqref{approxprob} against $\varphi$ to obtain
\begin{align*}
\Big|\intomega \partial_t\nep\varphi\Big|&=\Big|\intomega \Big(-\uep\cdot\!\nabla \nep+\Delta \nep-\nabla\!\cdot\big(\rhoep\fep(\nep)\nep\nabla \cep\big)+\kappa \nep-\mu \nep^2\Big)\varphi\Big|\\
&=\Big|\intomega \nep(\uep\cdot\nabla\varphi)-\intomega \nabla\nep\cdot\nabla\varphi+\intomega\rhoep\fep(\nep)\nep(\nabla \cep\cdot\nabla\varphi)+\kappa\intomega \nep\varphi-\mu\intomega \nep^2\varphi\Big|
\end{align*}
for all $t>0$ and $\epsi\in(0,1)$, where the boundary integrals again disappear due to the boundary conditions for $\nep$ and $\uep$ and the fact that $\rhoep=0$ on $\romega$. Here, we deduce from multiple applications of Hölder's inequality and the fact that for all $\epsi\in(0,1)$ we have $\rhoep f'(\nep)\leq 1$ on $\Omega\times(0,\infty)$ that
\begin{align*}
\Big|\intomega\partial_t\nep\varphi\Big|\leq  \big(\|\nep\|_{\Lo[2]}\|\uep\|_{\Lo[4]}&+\|\nabla\nep\|_{\Lo[\frac 43]}+\|\nep\|_{\Lo[2]}\|\nabla\cep\|_{\Lo[4]}\big)\|\nabla\varphi\|_{\Lo[4]}\\
&+\big(\kappa\|\nep\|_{\Lo[1]}+ \mu\|\nep\|_{\Lo[2]}^2\big)\|\varphi\|_{\Lo[\infty]}\quad\text{for all }t>0.
\end{align*}
Here, we make use of Lemma \ref{lem:u-stokes-bound} and multiple uses of Young's inequality to conclude that there is $C_1>0$ such that for all $\epsi\in(0,1)$
\begin{align*}
\Big|\intomega\partial_t\nep\varphi\Big|\leq C_1\big(\|\nep\|_{\Lo[2]}^2+\|\nabla\nep\|_{\Lo[\frac{4}{3}]}^\frac{4}{3}&+\|\nabla\cep\|_{\Lo[4]}^4+1\big)\|\nabla\varphi\|_{\Lo[4]}\\&+\big(\kappa\|\nep\|_{\Lo[1]}+\mu\|\nep\|_{\Lo[2]}^2\big)\|\varphi\|_{\Lo[\infty]}
\end{align*}
holds for all $t>0.$ Then, since $\|\varphi\|_{\LSpn{\infty}{(0,T);W_0^{1,4}(\Omega)}}\leq 1$ and $W_0^{1,4}(\Omega)\hookrightarrow\Lo[\infty]$, an integration over $(0,T)$ immediately entails \eqref{eq:n-time-regu} thanks to Lemmas \ref{lem:simple-bounds}, \ref{lem:grad-c-l4} and \ref{lem:grad-n-time}.

Similarly, fixing $\hat\varphi\in\LSpb{\infty}{(0,T);W_0^{1,4}(\Omega)}$ with $\|\hat\varphi\|_{\LSpn{\infty}{(0,T);W_0^{1,4}(\Omega)}}\leq1$ and testing the second equation of \eqref{approxprob} against $\hat\varphi$ we find $C_2>0$ satisfying
\begin{align*}
\intoT\Big|\intomega\partial_t\cep\hat\varphi\Big|&=\intoT\Big|\intomega\Big(-\uep\cdot\nabla\cep+\Delta\cep-\gep(\nep)\cep\Big)\hat\varphi\Big|\\
&\leq C_2\intoT\!\Big(\|\nabla\cep\|_{\Lo[4]}^4+\|\Delta\cep\|_{\Lo[2]}^2+\|\nep\|_{\Lo[2]}^2+1\Big)\|\hat\varphi\|_{\Lo[2]}
\end{align*}
for all $\epsi\in(0,1)$, where we again made use of Lemmas \ref{lem:u-stokes-bound} and \ref{lem:simple-bounds}, Young's inequality and the fact that $f(s)\leq s$ for all $s\geq0$, so that \eqref{eq:c-time-regu} is an evident consequence of the spatio-temporal bounds provided by Lemmas \ref{lem:grad-c-l4} and \ref{lem:simple-bounds}.

Finally, for any fixed $\psi\in \LSpb{2}{(0,T);W_{0,\sigma}^{1,2}(\Omega)}$ with $\|\psi\|_{\LSpn{2}{(0,T);W_{0,\sigma}^{1,2}(\Omega)}}\leq1$, we multiply the third equation in \ref{approxprob} by $\psi$ and integrate the resulting equation to derive that
\begin{align*}
&\intoT\Big|\intomega \partial_t\uep\cdot\psi\Big|^2=\int_0^T\Big|\intomega\nabla \uep\cdot\nabla\psi+\intomega \nep(\nabla\phi\cdot\psi)\Big|^2\\
\leq\ &C_3\intoT \|\nabla\uep\|_{\Lo[2]}^2\|\nabla\psi\|_{\Lo[2]}+C_3\intoT\|\nep\|_{\Lo[2]}^2\|\nabla\phi\|_{\Lo[\infty]}^2\|\psi\|^2_{\Lo[2]}\leq C_4(T)
\end{align*}
in light of Lemma \ref{lem:u-stokes-bound}, \eqref{phi-def} and Lemma \ref{lem:simple-bounds}, and from which we conclude \eqref{eq:u-time-regu}.

\end{bew}

\setcounter{equation}{0} 
\section{Existence of a limit solution. The proof of Theorem \ref{theo:1}}\label{sec6:limit}
Collecting the uniform bounds presented in Sections \ref{sec2:weak-sol}--\ref{sec5:time-reg}, we can now construct a limit object which satisfies all the regularity requirements present in Definition \ref{def:sol}.

\begin{proposition}
\label{prop:convergence}
There exist a sequence $(\epsi_j)_{j\in \N}\subset(0,1)$ with $\epsi_j\searrow0$ as $j\to\infty$ and functions
	\begin{align*}
n&\in \LSploc{2}{\bomega\times[0,\infty)}\quad\text{with}\quad\nabla n\in\LSploc{\frac43}{\bomega\times[0,\infty);\R^\dimN},\\
c&\in \LSp{\infty}{\Omega\times(0,\infty)}\quad\text{with}\quad c-\cs\in\LSplocb{4}{[0,\infty);W_0^{1,4}(\Omega)},\\
u&\in \LSplocb{2}{[0,\infty);W_{0,\sigma}^{1,2}(\Omega)}
\end{align*}
and such that the solutions $(\nep,\cep,\uep)$ of \eqref{approxprob} fulfill
\begin{alignat}{2}
\nep\to\ & n&&\text{in }\ \LSploc{p}{\bomega\times[0,\infty)}\text{ for any }p\in[1,2]\text{ and a.e. in}\ \Omega\times(0,\infty),\label{eq:conv-n}\\
\nep\wto\ &  n&&\text{in }\ \LSploc{2}{\bomega\times[0,\infty)},\label{eq:conv-weak-n}\\
\nabla\nep\wto\  &\nabla n&\ \; &\text{in }\ \LSploc{\frac43}{\bomega\times[0,\infty);\R^\dimN},\label{eq:conv-nab-n}\\
\rhoep\fep(\nep)\nep\to\ &n&\ &\text{in }\LSploc{p}{\bomega\times[0,\infty)}\text{ for any }p\in[1,2),\label{eq:conv-fep_n}\\
\gep(\nep)\to\ &n&\ &\text{in }\LSploc{p}{\bomega\times[0,\infty)}\text{ for any }p\in[1,2),\label{eq:conv-gep}\\
\cep\ \to\ & c&&\text{in }\ \LSploc{q}{\bomega\times[0,\infty)}\text{ for any }q\in[1,\infty)\text{ and a.e. in}\ \Omega\times(0,\infty),\label{eq:conv-c}\\
\cep\wsto\ & c&&\text{in }\ \LSp{\infty}{\Omega\times(0,\infty)},\label{eq:conv-weaks-c}\\
\nabla\cep \wto\ & \nabla c&&\text{in }\ \LSploc{4}{\bomega\times[0,\infty);\R^\dimN},\label{eq:conv-nab-c}\\
\uep\to\ & u&&\text{in }\ \LSploc{r}{\bomega\times[0,\infty);\R^\dimN}\text{ for any }r\in[1,6)\text{ and a.e. in}\ \Omega\times(0,\infty),\label{eq:conv-u}\\
\uep\wsto\ & u&&\text{in }\ \LSp{\infty}{(0,\infty);\LSp{6}{\Omega;\R^\dimN}},\label{eq:conv-weaks-u}\\
\nabla \uep \wto\ &\nabla u &&\text{in }\ \LSploc{2}{\bomega\times[0,\infty);\R^{\dimN\times \dimN}},\label{eq:conv-nab-u}
\end{alignat}
as $\epsi=\epsi_j\searrow0$.
\end{proposition}

\begin{bew}
A combination of an Aubin--Lions type lemma (\cite[Corollary 8.4]{Sim87}) with the bounds presented in Lemmas \ref{lem:simple-bounds}, \ref{lem:grad-n-time} and \ref{lem:time-reg} ensures that
\begin{align*}
\left\{\nep\right\}_{\epsi\in(0,1)}\quad\text{is relatively compact in }\LSploc{\frac43}{\bomega\times[0,\infty)}
\end{align*}
and that hence we can find a subsequence $(\epsi_j)_{j\in\N}$ with $\epsi_j\searrow0$ as $j\to\infty$ such that $\nep\to n$ in $\LSploc{\frac43}{\bomega\times[0,\infty)}$ and a.e. in $\Omega\times(0,\infty)$. Additionally, the spatio-temporal bounds of Lemmas \ref{lem:simple-bounds} and \ref{lem:grad-n-time} also allow us to conclude \eqref{eq:conv-weak-n} and \eqref{eq:conv-nab-n}, respectively, along a subsequence (still denoted by $(\epsi_j)_{j\in\N}$). Moreover, noting that $\Theta:[0,\infty)\to[0,\infty]$, $\Theta(s):=s\ln(s+1)$ is an increasing and convex function with $$\lim_{s\to\infty}\frac{\Theta(s)}{s}=\infty\quad\text{and}\quad\Theta(s)\leq 2s\ln(s^\frac{1}{2})+1\quad\text{for all }s\geq0,$$ we observe that according to Lemma \ref{lem:eq-n} for any $T>0$ there is $C>0$ satisfying
\begin{align*}
\intoTomega\Theta(\nep^2)=\intoTomega\nep^2\ln(\nep^2+1)\leq 2\intoTomega\nep^2\ln\nep+|\Omega|T\leq C(T),
\end{align*} 
which, by a result of de la Vallée--Poussin (e.g. \cite[II.T22]{dellacherieProbabilitiesPotential1978}), entails that $\left\{\nep^2\right\}_{\epsi\in(0,1)}$ is equi-integrable. Thus, a combination of the equi-integrability with the a.e. convergence of $\nep$ and Vitali's theorem yields \eqref{eq:conv-n} along a further subsequence. Then, since $|\rhoep\fep(\nep)|\leq 1$ in $\Omega\times(0,\infty)$ for all $\epsi\in(0,1)$ and $\rhoep\fep(\nep)\to 1$ a.e. in $\Omega\times(0,\infty)$ as $\epsi=\epsi_j\searrow0$, we find from \eqref{eq:conv-n} and arguments akin to e.g. \cite[Lemma A.4]{win15_chemorot} that \eqref{eq:conv-fep_n} holds as well. Likewise, we may also conclude \eqref{eq:conv-gep} from \eqref{eq:conv-n}. Working along similar lines for the second and third components, we can draw on the bounds of Lemmas \ref{lem:simple-bounds}, \ref{lem:grad-c-l4} and \ref{lem:time-reg} to obtain an additional subsequence along which \eqref{eq:conv-c}, \eqref{eq:conv-weaks-c} and \eqref{eq:conv-nab-c} hold and iterating the arguments once more with the bounds of Lemmas \ref{lem:u-stokes-bound} and \ref{lem:time-reg} concerning $\uep$, finally, also \eqref{eq:conv-u}, \eqref{eq:conv-weaks-u} and \eqref{eq:conv-nab-u}. The claimed regularity properties of $(n,c,u)$ and $c-\cs$ are clearly a direct consequence of \eqref{eq:conv-weak-n}, \eqref{eq:conv-nab-n}, \eqref{eq:conv-weaks-c}, \eqref{eq:conv-nab-c}, \eqref{eq:conv-u}, \eqref{eq:conv-nab-u} and \eqref{eq:cstar-reg} and the fact that $\cep-\cs=0$ on $\romega\times[0,\infty)$.
\end{bew}

Finally, it remains to be checked that the limit objected provided by Proposition \ref{prop:convergence} indeed satisfies the integral identities \eqref{eq:def-eq-n}, \eqref{eq:def-eq-c} and \eqref{eq:def-eq-u} of Definition \ref{def:sol}. This, however, is a straightforward procedure, as the convergence properties of Proposition \ref{prop:convergence} already cover everything we need to pass to the limit in the corresponding equations of \eqref{approxprob}.

\begin{proof}[\textbf{Proof of Theorem \ref{theo:1}}:]
Since the regularity requirements imposed on a weak solution by Definition \ref{def:sol} are already covered by the properties obtained in Lemma \ref{prop:convergence}, we only have to verify that the limit objects obtained in said lemma also satisfy the integral identities \eqref{eq:def-eq-n}--\eqref{eq:def-eq-u}.
We pick $\varphi\in C_0^\infty\big(\bomega\times[0,\infty)\big)$, $\hat{\varphi}\in C_0^\infty\big(\Omega\times[0,\infty)\big)$ and $\psi\in C_{0}^\infty\big(\bomega\times[0,\infty);\R^\dimN\big)$ with $\nabla\cdot\psi\equiv0$ and then fix $T>0$ such that $\varphi,\hat{\varphi},\psi\equiv0$ in $\Omega\times(T,\infty)$. Now, we test the first equation of \eqref{approxprob} against $\varphi$ and integrate by parts to obtain
\begin{align}\label{eq:theo1-eq1}
-\intoTomega \nep\varphi_t-\intomega n_0\varphi(\cdot,0)=-\intoTomega\nabla \nep\cdot\nabla\varphi&+\intinfomega \rhoep\fep(\nep)\nep (\nabla \cep\cdot\nabla\varphi)\\&+\kappa\intoTomega \nep\varphi-\mu\intoTomega \nep^2\varphi+\intoTomega \nep (\uep\cdot\nabla\varphi)\nonumber
\end{align}
for all $\epsi\in(0,1)$, where we made use of the fact that $\uep$ is solenoidal. According to \eqref{eq:conv-n} and \eqref{eq:conv-nab-n}, we immediately find that
\[
\intoTomega \nep\varphi_t\to \intoTomega n \varphi_t , \quad \intoTomega\nabla \nep\cdot\nabla\varphi\to \intoTomega\nabla n\cdot\nabla\varphi   
\]
and
\[
\intoTomega \nep\varphi\to \intoTomega n\varphi   \quad \mathrm{as}\,\,\epsi=\epsi_j\searrow 0.
\]
Drawing on \eqref{eq:conv-n}, \eqref{eq:conv-fep_n} combined with \eqref{eq:conv-nab-c} and \eqref{eq:conv-gep} together with \eqref{eq:conv-u} we also conclude

\[
\intoTomega \nep^2\varphi\to \intoTomega n^2\varphi, \quad\quad \intoTomega \rhoep\fep(\nep)\nep (\nabla \cep\cdot\nabla\varphi)\to    \intoTomega n (\nabla c\cdot\nabla\varphi)
\]
and
\[
\intoTomega \nep (\uep\cdot\nabla\varphi)   \to        \intoTomega n (u\cdot\nabla\varphi)   \quad \mathrm{as}\,\,\epsi=\epsi_j\searrow 0
\]
so that passing to the limit in \eqref{eq:theo1-eq1} immediately entails \eqref{eq:def-eq-n} since $T>0$ was chosen such that $\varphi\equiv 0$ in $\Omega\times(T,\infty)$.

Next, multiplying the second equation of \eqref{approxprob} by $\hat{\varphi}$ and integrating, we have
\begin{align*}
-\intoTomega \cep\hat{\varphi}_t-\intomega c_0\hat{\varphi}(\cdot,0)=-\intoTomega\nabla \cep\cdot\nabla\hat{\varphi}-\intoTomega \gep(\nep)\cep\hat{\varphi}+\intoTomega \cep (\uep\cdot\nabla\hat{\varphi})
\end{align*}
for all $\epsi\in(0,1)$. Therefore, utilizing \eqref{eq:conv-c}, \eqref{eq:conv-nab-c}, \eqref{eq:conv-gep} and \eqref{eq:conv-u}, we may also pass to the limit in this equation and obtain \eqref{eq:def-eq-c}.
Finally, testing the third equation in \eqref{approxprob} against $\psi$ and integrating by parts yields
\begin{align*}
-\intoTomega \uep\cdot\psi_t-\intomega u_0\cdot\psi(\cdot,0)=-\intoTomega\nabla \uep\cdot\nabla\psi+\intoTomega \nep(\nabla\phi\cdot\psi)\quad\text{for all }\epsi\in(0,1),
\end{align*}
where \eqref{eq:conv-u}, \eqref{eq:conv-nab-u} and \eqref{eq:conv-n} imply that we may, once more, pass to the limit and obtain \eqref{eq:def-eq-u}, concluding the proof.
\end{proof}
{\small

\section*{Acknowledgements}
Tobias Black acknowledges support of the {\em Deutsche Forschungsgemeinschaft} in the context of the project
  {\em Emergence of structures and advantages in cross-diffusion systems (project no. 411007140)}. Chunyan Wu was sup-
ported by the {\em CSC (China Scholarship Council)} and the {\em Applied Fundamental Research Program of Sichuan
Province (no. 2020YJ0264)} during her stay at the University of Paderborn.

\footnotesize{
\setlength{\bibsep}{2pt plus 0.5ex}

\begin{thebibliography}{46}
\providecommand{\natexlab}[1]{#1}

\bibitem[Bellomo et~al.(2015)Bellomo, Bellouquid, Tao, and
  Winkler]{BBWT15}{https://doi.org/10.1142/S021820251550044X}
N.~Bellomo, A.~Bellouquid, Y.~Tao, and M.~Winkler.
\newblock Toward a mathematical theory of {K}eller-{S}egel models of pattern
  formation in biological tissues.
\newblock \emph{Math. Models Methods Appl. Sci.}, 25\penalty0 (9):\penalty0
  1663--1763, 2015.

\bibitem[Braukhoff(2017)]{braukhoffGlobalWeakSolution2017}{https://doi.org/10.1016/j.anihpc.2016.08.003}
M.~Braukhoff.
\newblock Global (weak) solution of the chemotaxis-{Navier}–{Stokes}
  equations with non-homogeneous boundary conditions and logistic growth.
\newblock \emph{Ann. Inst. H. Poincaré Anal. Non Linéaire}, 34\penalty0
  (4):\penalty0 1013--1039, 2017.

\bibitem[Braukhoff and
  Lankeit(2019)]{braukhoffStationarySolutionsChemotaxisconsumption2019}{https://doi.org/10.1142/S0218202519500398}
M.~Braukhoff and J.~Lankeit.
\newblock Stationary solutions to a chemotaxis-consumption model with realistic
  boundary conditions for the oxygen.
\newblock \emph{Math. Models Methods Appl. Sci.}, 29\penalty0 (11):\penalty0
  2033--2062, 2019.

\bibitem[Cao and
  Lankeit(2016)]{caolan16_smalldatasol3dnavstokes}{https://doi.org/10.1007/s00526-016-1027-2}
X.~Cao and J.~Lankeit.
\newblock Global classical small-data solutions for a three-dimensional
  chemotaxis {Navier}-{Stokes} system involving matrix-valued sensitivities.
\newblock \emph{Calc. Var. Partial Differential Equations}, 55\penalty0
  (4):\penalty0 Paper No. 107, 2016.

\bibitem[Chertock et~al.(2012)Chertock, Fellner, Kurganov, Lorz, and
  Markowich]{chertockSinkingMergingStationary2012}{https://doi.org/10.1017/jfm.2011.534}
A.~Chertock, K.~Fellner, A.~Kurganov, A.~Lorz, and P.~A. Markowich.
\newblock Sinking, merging and stationary plumes in a coupled chemotaxis-fluid
  model: a high-resolution numerical approach.
\newblock \emph{Journal of Fluid Mechanics}, 694:\penalty0 155--190, 2012.

\bibitem[Dellacherie and Meyer(1978)]{dellacherieProbabilitiesPotential1978}{}
C.~Dellacherie and P.~A. Meyer.
\newblock \emph{Probabilities and potential}.
\newblock Number~29 in North-{Holland} mathematics studies. Hermann, Paris,
  1978.

\bibitem[Dombrowski et~al.(2004)Dombrowski, Cisneros, Chatkaew, Goldstein, and
  Kessler]{Dombrowski04}{https://doi.org/10.1103/PhysRevLett.93.098103}
C.~Dombrowski, L.~Cisneros, S.~Chatkaew, R.~E. Goldstein, and J.~O. Kessler.
\newblock Self-{Concentration} and {Large}-{Scale} {Coherence} in {Bacterial}
  {Dynamics}.
\newblock \emph{Physical Review Letters}, 93\penalty0 (9), 2004.

\bibitem[Duan et~al.(2010)Duan, Lorz, and
  Markowich]{Duan2010}{https://doi.org/10.1080/03605302.2010.497199}
R.~Duan, A.~Lorz, and P.~Markowich.
\newblock Global {Solutions} to the {Coupled} {Chemotaxis}-{Fluid} {Equations}.
\newblock \emph{Comm. Partial Differential Equations}, 35\penalty0
  (9):\penalty0 1635--1673, 2010.

\bibitem[Evans(2010)]{evans}{https://doi.org/10.1090/gsm/019}
L.~C. Evans.
\newblock \emph{Partial differential equations}, Volume~19 of \emph{Graduate
  {Studies} in {Mathematics}}.
\newblock American Mathematical Society, Providence, RI, second edition, 2010.

\bibitem[Fuest et~al.(2021)Fuest, Lankeit, and
  Mizukami]{fuestLongtermBehaviourParabolic2021}{https://doi.org/10.1016/j.jde.2020.08.021}
M.~Fuest, J.~Lankeit, and M.~Mizukami.
\newblock Long-term behaviour in a parabolic–elliptic
  chemotaxis–consumption model.
\newblock \emph{J. Differential Equations}, 271:\penalty0 254--279, 2021.

\bibitem[Fujikawa and
  Matsushita(1989)]{fujikawaFractalGrowthBacillus1989}{https://doi.org/10.1143/JPSJ.58.3875}
H.~Fujikawa and M.~Matsushita.
\newblock Fractal {Growth} of \textit{{Bacillus} subtilis} on {Agar} {Plates}.
\newblock \emph{J. Phys. Soc. Japan}, 58\penalty0 (11):\penalty0 3875--3878,
  1989.

\bibitem[Giga(1986)]{gig86}{https://doi.org/10.1016/0022-0396(86)90096-3}
Y.~Giga.
\newblock Solutions for semilinear parabolic equations in {$L^p$} and
  regularity of weak solutions of the {Navier}-{Stokes} system.
\newblock \emph{J. Differential Equations}, 62\penalty0 (2), 1986.

\bibitem[Henry(1981)]{hen81}{https://doi.org/10.1007/BFb0089647}
D.~Henry.
\newblock \emph{Geometric {Theory} of {Semilinear} {Parabolic} {Equations}},
  Volume 840 of \emph{Lecture {Notes} in {Mathematics}}.
\newblock Springer Berlin Heidelberg, 1981.

\bibitem[Horstmann(2003)]{Ho03}{}
D.~Horstmann.
\newblock From 1970 until present: the {Keller}-{Segel} model in chemotaxis and
  its consequences. {I}.
\newblock \emph{Jahresber. Deutsch. Math.-Verein.}, 105\penalty0 (3):\penalty0
  103--165, 2003.

\bibitem[Keller and
  Segel(1970)]{KS70}{https://doi.org/10.1016/0022-5193(70)90092-5}
E.~F. Keller and L.~A. Segel.
\newblock Initiation of slime mold aggregation viewed as an instability.
\newblock \emph{J. Theor. Biol.}, 26\penalty0 (3):\penalty0 399--415, 1970.
\newblock tex.fjournal: Journal of Theoretical Biology.

\bibitem[Knosalla(2017)]{knosallaGlobalSolutionsAerotaxis2017}{https://doi.org/10.4064/am2301-2-2017}
P.~Knosalla.
\newblock Global solutions of aerotaxis equations.
\newblock \emph{Applicationes Mathematicae}, 44\penalty0 (1):\penalty0
  135--148, 2017.

\bibitem[Ladyženskaja et~al.(1968)Ladyženskaja, Solonnikov, and
  Ural'ceva]{LSU}{https://doi.org/10.1090/mmono/023}
O.~A. Ladyženskaja, V.~A. Solonnikov, and N.~N. Ural'ceva.
\newblock \emph{Linear and quasilinear equations of parabolic type}.
\newblock Translations of mathematical monographs. American Mathematical
  Society, 1968.

\bibitem[Lankeit and
  Lankeit(2019)]{LL18-ClassSolLogCTSingSens}{https://doi.org/10.1016/j.nonrwa.2018.09.012}
E.~Lankeit and J.~Lankeit.
\newblock Classical solutions to a logistic chemotaxis model with singular
  sensitivity and signal absorption.
\newblock \emph{Nonlinear Anal. Real World Appl.}, 46:\penalty0 421--445, 2019.

\bibitem[Lankeit(2015)]{Lan15-threshold}{https://doi.org/10.3934/dcdsb.2015.20.1499}
J.~Lankeit.
\newblock Chemotaxis can prevent thresholds on population density.
\newblock \emph{Discrete Contin. Dyn. Syst. Ser. B}, 20\penalty0 (5):\penalty0
  1499--1527, 2015.

\bibitem[Lankeit(2016)]{Lan-Longterm_M3AS16}{https://doi.org/10.1142/S021820251640008X}
J.~Lankeit.
\newblock Long-term behaviour in a chemotaxis-fluid system with logistic
  source.
\newblock \emph{Math. Models Methods Appl. Sci.}, 26\penalty0 (11):\penalty0
  2071--2109, 2016.

\bibitem[Lankeit and
  Winkler(2019)]{LanWin_JDMV_19}{https://doi.org/10.1365/s13291-019-00210-z}
J.~Lankeit and M.~Winkler.
\newblock Facing low regularity in chemotaxis systems.
\newblock \emph{Jahresber. Deutsch. Math.-Verein.}, 2019.

\bibitem[Lee and
  Kim(2015)]{leeNumericalInvestigationFalling2015}{https://doi.org/10.1016/j.euromechflu.2015.03.002}
H.~G. Lee and J.~Kim.
\newblock Numerical investigation of falling bacterial plumes caused by
  bioconvection in a three-dimensional chamber.
\newblock \emph{Eur. J. of Mech. B Fluids}, 52:\penalty0 120--130, 2015.

\bibitem[Lorz(2010)]{lorz-M3AS10}{https://doi.org/10.1142/S0218202510004507}
A.~Lorz.
\newblock Coupled {Chemotaxis} {Fluid} {Model}.
\newblock \emph{Math. Models Methods Appl. Sci.}, 20\penalty0 (06):\penalty0
  987--1004, 2010.

\bibitem[Matsushita and
  Fujikawa(1990)]{matsushitaDiffusionlimitedGrowthBacterial1990}{https://doi.org/10.1016/0378-4371(90)90402-E}
M.~Matsushita and H.~Fujikawa.
\newblock Diffusion-limited growth in bacterial colony formation.
\newblock \emph{Physica A: Statistical Mechanics and its Applications},
  168\penalty0 (1):\penalty0 498--506, 1990.

\bibitem[Peng and
  Xiang(2018)]{pengGlobalSolutionsCoupled2018}{https://doi.org/10.1142/S0218202518500239}
Y.~Peng and Z.~Xiang.
\newblock Global solutions to the coupled chemotaxis-fluids system in a {3D}
  unbounded domain with boundary.
\newblock \emph{Math. Models Methods Appl. Sci.}, 28\penalty0 (05):\penalty0
  869--920, 2018.

\bibitem[Peng and
  Xiang(2019)]{pengGlobalExistenceConvergence2019}{https://doi.org/10.1016/j.jde.2019.02.007}
Y.~Peng and Z.~Xiang.
\newblock Global existence and convergence rates to a chemotaxis-fluids system
  with mixed boundary conditions.
\newblock \emph{J. Differential Equations}, 267\penalty0 (2):\penalty0
  1277--1321, 2019.

\bibitem[Porzio and
  Vespri(1993)]{PorzVesp93}{https://doi.org/10.1006/jdeq.1993.1045}
M.~M. Porzio and V.~Vespri.
\newblock Hölder estimates for local solutions of some doubly nonlinear
  degenerate parabolic equations.
\newblock \emph{J. Differential Equations}, 103\penalty0 (1):\penalty0
  146--178, 1993.

\bibitem[Quittner and
  Souplet(2007)]{QS07}{https://doi.org/10.1007/978-3-030-18222-9}
P.~Quittner and P.~Souplet.
\newblock \emph{Superlinear parabolic problems}.
\newblock Birkhäuser {Advanced} {Texts}: {Basler} {Lehrbücher}. Birkhäuser
  Verlag, Basel, 2007.

\bibitem[Simon(1987)]{Sim87}{https://doi.org/10.1007/BF01762360}
J.~Simon.
\newblock Compact sets in the space {$L^p(0,T;B)$}.
\newblock \emph{Ann. Mat. Pura Appl. (4)}, 146:\penalty0 65--96, 1987.

\bibitem[Sohr(2001)]{sohr}{https://doi.org/10.1007/978-3-0348-8255-2}
H.~Sohr.
\newblock \emph{The {Navier}-{Stokes} equations}.
\newblock Birkhäuser {Advanced} {Texts}: {Basler} {Lehrbücher}. Birkhäuser
  Verlag, Basel, 2001.

\bibitem[Solonnikov(2007)]{Solonnikov2007}{https://doi.org/10.1090/trans2/220/08}
V.~A. Solonnikov.
\newblock Schauder estimates for the evolutionary generalized {S}tokes problem.
\newblock In \emph{Nonlinear equations and spectral theory}, Volume 220 of
  \emph{Amer. Math. Soc. Transl. Ser. 2}, pp. 165--200. Amer. Math. Soc.,
  Providence, RI, 2007.

\bibitem[Tao and
  Winkler(2012)]{TaoWin-quasilinear_JDE12}{https://doi.org/10.1016/j.jde.2011.08.019}
Y.~Tao and M.~Winkler.
\newblock Boundedness in a quasilinear parabolic-parabolic {Keller}-{Segel}
  system with subcritical sensitivity.
\newblock \emph{J. Differential Equations}, 252\penalty0 (1), 2012.

\bibitem[Tuval et~al.(2005)Tuval, Cisneros, Dombrowski, Wolgemuth, Kessler, and
  Goldstein]{tuval2005bacterial}{https://doi.org/10.1073/pnas.0406724102}
I.~Tuval, L.~Cisneros, C.~Dombrowski, C.~W. Wolgemuth, J.~O. Kessler, and R.~E.
  Goldstein.
\newblock Bacterial swimming and oxygen transport near contact lines.
\newblock \emph{Proc. Natl. Acad. Sci. U.S.A.}, 102\penalty0 (7):\penalty0
  2277--2282, 2005.

\bibitem[Vorotnikov(2014)]{Vorotnikov-WeakSol-CMS14}{https://doi.org/10.4310/CMS.2014.v12.n3.a8}
D.~Vorotnikov.
\newblock Weak solutions for a bioconvection model related to {Bacillus}
  subtilis.
\newblock \emph{Commun. Math. Sci.}, 12\penalty0 (3):\penalty0 545--563, 2014.

\bibitem[Wang et~al.(2020{\natexlab{a}})Wang, Winkler, and
  Xiang]{WWX2019-global_mass-preserving_dirichlet_signal}{}
Y.~Wang, M.~Winkler, and Z.~Xiang.
\newblock Global mass-preserving solutions to a chemotaxis-fluid model
  involving {D}irichlet boundary conditions for the signal.
\newblock 2020{\natexlab{a}}.
\newblock Preprint.

\bibitem[Wang et~al.(2020{\natexlab{b}})Wang, Winkler, and
  Xiang]{WWX2019-local_energy_estimates_prescribed_signal}{}
Y.~Wang, M.~Winkler, and Z.~Xiang.
\newblock Local energy estimates and global solvability in a three-dimensional
  chemotaxis-fluid system with prescribed signal on the boundary.
\newblock 2020{\natexlab{b}}.
\newblock Preprint.

\bibitem[Winkler(2010)]{win10jde}{https://doi.org/10.1016/j.jde.2010.02.008}
M.~Winkler.
\newblock Aggregation vs. global diffusive behavior in the higher-dimensional
  {Keller}-{Segel} model.
\newblock \emph{J. Differential Equations}, 248\penalty0 (12):\penalty0
  2889--2905, 2010.

\bibitem[Winkler(2012)]{win_fluid_CPDE12}{https://doi.org/10.1080/03605302.2011.591865}
M.~Winkler.
\newblock Global large-data solutions in a chemotaxis-({Navier}-){Stokes}
  system modeling cellular swimming in fluid drops.
\newblock \emph{Comm. Partial Differential Equations}, 37\penalty0
  (2):\penalty0 319--351, 2012.

\bibitem[Winkler(2014)]{win-stab2d-ArchRatMechAna12}{https://doi.org/10.1007/s00205-013-0678-9}
M.~Winkler.
\newblock Stabilization in a two-dimensional chemotaxis-{Navier}-{Stokes}
  system.
\newblock \emph{Arch. Ration. Mech. Anal.}, 211\penalty0 (2):\penalty0
  455--487, 2014.

\bibitem[Winkler(2015{\natexlab{a}})]{Win-ct_fluid_3d-CPDE15}{https://doi.org/10.1007/s00526-015-0922-2}
M.~Winkler.
\newblock Boundedness and large time behavior in a three-dimensional
  chemotaxis-{Stokes} system with nonlinear diffusion and general sensitivity.
\newblock \emph{Calc. Var. Partial Differential Equations}, 54\penalty0
  (4):\penalty0 3789--3828, 2015{\natexlab{a}}.

\bibitem[Winkler(2015{\natexlab{b}})]{win15_chemorot}{https://doi.org/10.1137/140979708}
M.~Winkler.
\newblock Large-data global generalized solutions in a chemotaxis system with
  tensor-valued sensitivities.
\newblock \emph{SIAM J. Math. Anal.}, 47\penalty0 (4):\penalty0 3092--3115,
  2015{\natexlab{b}}.

\bibitem[Winkler(2016)]{win_globweak3d-AHPN16}{https://doi.org/10.1016/j.anihpc.2015.05.002}
M.~Winkler.
\newblock Global weak solutions in a three-dimensional
  chemotaxis--{Navier}-{Stokes} system.
\newblock \emph{Ann. Inst. H. Poincaré Anal. Non Linéaire}, 33\penalty0
  (5):\penalty0 1329--1352, 2016.

\bibitem[Winkler(2017)]{win_chemonavstokesfinal_TransAm17}{https://doi.org/10.1090/tran/6733}
M.~Winkler.
\newblock How far do chemotaxis-driven forces influence regularity in the
  {Navier}-{Stokes} system?
\newblock \emph{Trans. Amer. Math. Soc.}, 369\penalty0 (5):\penalty0
  3067--3125, 2017.

\bibitem[Winkler(2020)]{win20-leray-structure-ct-fluid}{}
M.~Winkler.
\newblock Does {Leray's} structure theorem withstand buoyancy-driven
  chemotaxis-fluid interaction?
\newblock 2020.
\newblock To appear in J. Eur. Math. Soc. (JEMS).

\bibitem[Woodward et~al.(1995)Woodward, Tyson, Myerscough, Murray, Budrene, and
  Berg]{woodwardSpatiotemporalPatternsGenerated1995}{https://doi.org/10.1016/S0006-3495(95)80400-5}
D.~Woodward, R.~Tyson, M.~Myerscough, J.~Murray, E.~Budrene, and H.~Berg.
\newblock Spatio-temporal patterns generated by {Salmonella} typhimurium.
\newblock \emph{Biophysical Journal}, 68\penalty0 (5):\penalty0 2181--2189,
  1995.

\bibitem[Wu and
  Xiang(2020)]{wuAsymptoticDynamicsChemotaxisNavier2020}{https://doi.org/10.1142/S0218202520500244}
C.~Wu and Z.~Xiang.
\newblock Asymptotic dynamics on a chemotaxis-{Navier}–{Stokes} system with
  nonlinear diffusion and inhomogeneous boundary conditions.
\newblock \emph{Math. Models Methods Appl. Sci.}, 30\penalty0 (07):\penalty0
  1325--1374, 2020.

\end{thebibliography}

}
\end{document}